\documentclass{elsart}
\usepackage{amssymb}

{\theoremstyle{plain}%

{\theoremstyle{remark}

}
{\theoremstyle{definition}

\newtheorem{example}{Example}
}

\def\be{\begin{equation}} \def\ee{\end{equation}}
\def\bena{\begin{eqnarray*}} \def\ena{\end{eqnarray*}}
 
\def\l|{\left|} \def\r|{\right|}

\begin{document}

\begin{frontmatter}

\title{Convoluted $C$-cosine functions and semigroups. Relations with ultradistribution and hyperfunction sines}

\author{M. Kosti\' c and  S.\ Pilipovi\'c
}\thanks{This research was supported by MNZ\v ZS of Serbia, project
no. 144016}

\address{Faculty of Technical Sciences,
Department of Mathematics and Informatics,\\ University of Novi Sad,
21000 Novi Sad, Serbia\\
pilipovic@im.ns.ac.yu\\
Tel. + 381 21 485 2850, Fax: + 381 21  6350  458}



\begin{abstract}
Convoluted $C$-cosine functions and semigroups in a Banach space
setting extending the classes of fractionally integrated
$C$-cosine functions and semigroups are systematically analyzed.
Structural properties of such operator families are obtained.
Relations between convoluted $C$-cosine functions and analytic
convoluted $C$-semigroups, introduced and investigated in this
paper are given through the convoluted version of the abstract
Weierstrass formula which is also proved in the paper.
Ultradistribution and hyperfunction sines are connected with
analytic convoluted semigroups and ultradistribution semigroups.
Several examples of operators generating convoluted cosine
functions, (analytic) convoluted semigroups as well as
hyperfunction and ultradistribution sines illustrate the abstract
approach of the authors. As an application, it is proved that the
polyharmonic operator $(-\Delta)^{2^{n}},$ $n\in {\mathbb N},$
acting on $L^{2}[0,\pi]$ with appropriate boundary conditions,
generates an exponentially bounded $K_{n}$-convoluted cosine
function, and consequently, an exponentially bounded analytic
$K_{n+1}$-convoluted semigroup of angle $\frac{\pi}{2},$ for
suitable exponentially bounded kernels $K_{n} $ and $K_{n+1}.$
\end{abstract}



\end{frontmatter}

\section{Introduction and preliminaries}\label{intro}

\noindent We study a class of convoluted $C$-cosine functions
extending the class of $\alpha$-times integrated $C$-cosine
functions, $\alpha>0$ and continue our researches in \cite{koj}-
\cite{ms} where we investigate different kinds of convoluted
operator type families and their relations with (tempered)
ultradistribution semigroups and (Fourier) hyperfunction semigroups.

Local convoluted $C$-semigroups were introduced and studied in the
papers of I. Cior\u anescu and G. Lumer \cite{c51}-\cite{c73} who
related them to ultradistribution semigroups, in the particular case
$C=I$. We refer to \cite{b42}, \cite{cha}, \cite{ci1}, \cite{cizi},
\cite{er}, \cite{err}, \cite{keya}, \cite{k92}, \cite{mp1},
\cite{ms}, \cite{ku113} and \cite{me152} for further information
concerning ultradistribution semigroups. We analyze in this paper
ultradistribution and hyperfunction sines continuing the researches
of H. Komatsu \cite{k92} and P. C. Kunstmann \cite{ku113}.

A class of exponentially bounded convoluted semigroups is
introduced and studied  in \cite{valent} via the operator valued
Laplace transform while global convoluted semigroups which are not
necessarily exponentially bounded have been recently analyzed in
\cite{koj} and \cite{mn}. We also refer to investigations of B.
B\" aumer, G. Lumer and F. Neubrander, \cite{b40} and \cite{luno},
for the use of the asymptotic Laplace transform in the theory of
convoluted semigroups, as well as to the paper \cite{miler} of C.
M\"  uler for the approximations of local convoluted semigroups.
In this paper, we further study convoluted $C$-cosine functions
introduced in \cite{koj} and obtain several generalizations of
results known for integrated $C$-cosine functions (cf. \cite{a23},
\cite{kikici}, \cite{keya}, \cite{kikic}, \cite {ko99},
\cite{mp1}, \cite{miki}, \cite{sh1}, \cite{w366}, \cite{x263},
\cite{z350} and \cite{z360}). We analyze in Section 2
$K$-convoluted $C$-cosine functions by a trustworthy passing to
the theory of $K$-convoluted $C$-semigroups on product spaces and
we compare corresponding integral generators of such operator
families. Such an approach enables one to obtain several
properties of subgenerators of convoluted $C$-cosine functions. We
also focus our attention to the case $C=I$ and continue the
analysis of P. C. Kunstmann \cite{ku101} concerning stationary
dense operators in Banach spaces. We prove that every generator
$A$ of a (local) $\alpha$-times integrated cosine function is
stationary dense and satisfies $n(A)\leq \lfloor \frac{\lceil
\alpha \rceil+1}{2}\rfloor .$ It seems to be an open problem to
improve this inequality; nevertheless, the concept of
stationarity, whose application in the problems of maximal
regularity of abstract Cauchy problems is not clearly
understandable, makes a difference between integrated operator
type families and convoluted operator type families. We generalize
in Section 3 results of \cite{koj} which are related to the
Laplace transform of exponentially bounded $K$-convoluted
$C$-cosine functions in order to use them in the later analysis of
the polyharmonic operator $\Delta^{2^{n}}$ on $L^{2}[0,\pi]$.

Our main results are given in Sections 4, 5 and 6. In Section 4,
we obtain the Hille-Yosida type theorems for generators of
analytic convoluted $C$-semigroups introduced in this paper (see
also \cite{mn}) and prove the convoluted version of the abstract
Weierstrass formula connecting analytic convoluted $C$-semigroups
and convoluted $C$-cosine functions. We relate in Section 5
ultradistribution and hyperfunction sines to analytic convoluted
semigroups and note, in ultradistribution case, some differences
between Beurling and Roumieu-type ultradistribution sines. Theorem
15 connects  ultradistribution sines of $(M_{p})-$class, resp.,
$\{M_{p}\}-$class, with ultradistribution semigroups of
$(M_{p}^{2})-$class, resp., $\{M_{p}^{2}\}-$class. In the rest of
Section 5, we analyze relations between (local) integrated cosine
functions as well as convoluted cosine functions with
ultradistribution semigroups. Such results were firstly obtained
by V. Keyantuo in ~\cite[Theorem 3.1]{keya} and this theorem has
been recently generalized and analyzed  in ~\cite[Theorem 4.3,
Example 4.4]{mp1}. Our results can be used in the analysis of
abstract Cauchy problems in the framework of various vector-valued
generalized function spaces.

We discuss in examples of Section 6 the polyharmonic operators
acting on $L^{2}[0,\pi]$ and point out, motivated by \cite{b40},
situations when the theory of convoluted cosine functions and
semigroups ($C=I$) cannot be used in the analysis of a wide class of
elliptic differential operators acting on $L^{p}-$type spaces (cf.
E. B. Davies \cite{dav1}, \cite{dav2}). In order to prove that the
polyharmonic operator $\Delta^{2^{n}}$ on $L^{2}[0,\pi]$ generates a
convoluted cosine function, we essentially use the fact that
$-\Delta^{2n}$ generates an analytic $C_{0}$-semigroup of angle
$\frac{\pi}{2}$ proved by J. A. Goldstein in \cite{goja}, see also
~\cite[Example 24.11]{l1}. Still, it is an open problem to
characterize polynomials of $-\Delta$ in the framework of the theory
of convoluted cosine functions and semigroups. We refer to
~\cite[Section VIII, XXIV]{l1} for the  application of entire
regularized groups in the analysis of such kind of problems.
Following R. Beals \cite{b41}, \cite{b42}, we construct an
illustrative example of an operator $A$ acting on the Hardy space
$H^{p}({{\mathbb C}_{+}}),$ $1\leq p <\infty$ which generates a
hyperfunction sine, but not an ultradistribution sine. Local
integrated semigroups generated by multiplication operators were
explicitly constructed by W. Arendt, O. El-Mennaoui and V. Keyantuo
in \cite{a22} (cf. \cite{ko99} for integrated cosine functions). We
construct convoluted cosine functions generated by multiplication
operators in Example 3 where we also discuss the maximal interval of
existence of a convoluted cosine function and present an example of
a global non-exponentially bounded convoluted cosine function.

In order to concentrate the exposition on our main results,
several structural properties of $K$-convoluted $C$-semigroups and
cosine functions are given in the Appendix, see \cite{koj} for
more details. Because of that, we do not analyze composition
properties, perturbations and approximation type results for
convoluted $C$-cosine functions as well as the corresponding
abstract Cauchy problems. These themes will be treated in a
separate paper.

{\sc Notation}. By  $E$ and $L(E)$ are denoted a complex Banach
space and Banach algebra of bounded linear operators on $E$. For a
closed linear operator $A$ on $E$,  $D(A)$, $\hbox{Kern(A)}$,
$R(A)$, $\rho (A)$ denote its domain, kernel, range and resolvent
set, respectively. Put $D_{\infty }(A) := \bigcap \limits_{n\in
{{\mathbb N}_{0}}} D(A^{n}).$ By $[D(A)]$ is denoted the Banach
space $D(A)$ endowed with the graph norm. In this paper, $C\in L(E)$
is an injective operator satisfying $CA\subset AC.$

We recall the basic facts from the Denjoy-Karleman-Komatsu theory of
ultradistributions although a great part of our results can be
transferred to the case of $\omega$-type ultradistributions. In the
sequel, $(M_p)_p$ is a sequence of positive numbers, $M_0=1$, such
that the following conditions are satisfied: \vspace{0.15cm}
\newline $ (M.1)\; \; M_p^2\leq M_{p+1} M_{p-1},\ p\in {\mathbb N},$
\newline $(M.2)\;\;M_{n} \leq AH^n\min_{p+q=n}M_pM_q,\ n\in
{\mathbb N},\;\mbox{for some }A,\ H>0,$
\newline $(M.3)'\;\; \sum_{p=1}^{\infty}\frac{M_{p-1}}{M_{p}}<\infty$.

If $(M_{p})$ is such a sequence, then as a matter of routine, one
can check that $(M^2_p)$ also satisfies $(M.1),$ $(M.2)$ and
$(M.3)'.$

If $s>1$ then the Gevrey sequences $(p!^s)_p,$ $(p^{ps})_p$ or
$(\Gamma (1+ps))_p$ satisfy the above conditions. The associated
function is defined by $M(\rho):=\sup_{p\in{\mathbb N}}\ln
\frac{\rho^p}{M_p}, \rho
>0;\ M(0):=0.$ If $\lambda \in {\mathbb C},$ then
$M(\lambda):=M(|\lambda|).$

We refer to \cite{k91} and \cite{k82} for the basic properties of
locally convex space-valued ultradifferentiable functions defined on
${\mathbb R}$ and corresponding ultradistributions of the Beurling,
resp., Roumieu type. The classes of Beurling, resp., Roumieu
ultradistributions with values in a Banach space $E$ are denoted by
${\mathcal D}^{'(M_{p})}(E)$, resp., ${\mathcal D}^{'\{M_{p}\}}(E)$
or simply ${\mathcal D}^{'(M_{p})},$ resp., ${\mathcal
D}^{'\{M_{p}\}}$ in the case $E={\mathbb R}.$ We denote by $\ast$
either $(M_p)$ or $\{M_p\}$. The similar terminology is used for the
spaces of Beurling and Roumieu type ultradifferentiable functions.
The space of all scalar-valued ultradistributions of $\ast$-class
with the support contained in $[0,\infty)$ is denoted by ${\mathcal
D}^{' *}_{0}$ (${\mathcal D}^{' *}_{0}(E)$ in the case of $E-$valued
ultradistributions).

The spaces of tempered ultradistributions of Beurling and Roumieu
type (cf. \cite{G}, \cite{ms} and \cite{pilip}) are defined as duals of\\
$ {\mathcal S}^{(M_p)}({\mathbb R}):=\mbox{proj}\lim_{k\rightarrow
\infty}{\mathcal S}^{M_p,k}({\mathbb R}),\mbox{ resp., }{\mathcal
S}^{\{M_p\}}({\mathbb R}):=\mbox{ind}\lim_{k\rightarrow 0}{\mathcal
S}^{M_p,k}({\mathbb R}), $ \\ where ${\mathcal S}^{M_p,k}({\mathbb
R}):=\{\phi\in C^\infty({\mathbb R}) :\ ||\phi||_k<\infty\}, $ $k>0$ and \\
$ ||\phi||_k:=$sup$\{\frac{k^{\alpha+\beta}}{M_\alpha
M_\beta}(1+|t|^2)^{\beta/2}|\phi^{(\alpha)}(t)| : \ t\in {\mathbb
R},\ \alpha,\ \beta \in {{\mathbb N}_0} \}. $ We refer to the book
of A. Kaneko \cite{kan} for the basic facts about hyperfunctions
and Fourier  hyperfunctions.

{\sc Terminology used in the paper}.\vspace{0.1cm}

1. $LT(\mathbb C)$ denotes the space of all Laplace transforms of
locally integrable, exponentially bounded functions.

2. If $\omega >0,$ put $\Pi _{\omega }:=\{z\in {\mathbb C} :
Rez>\omega ^{2}-\frac{(Imz)^{2}}{4\omega ^{2}}\}.$ Note,\\
$\Pi_{\omega}=\{z^{2}: \; z\in {\mathbb C},\ Rez>\omega\}.$

3. Let $\varepsilon>0$ and $C_{\varepsilon}>0.$ The next region was
introduced by S. \={O}uchi in \cite{o192}: $
\Omega_{\varepsilon,C_{\varepsilon}}:=\{\lambda \in {\mathbb C}:\;
Re\lambda \geq \varepsilon |\lambda|+C_{\varepsilon} \}.$ We will
use the notation\\ $\Omega^{2}_{\varepsilon, C_{\varepsilon}
}:=\{\lambda^{2}: \lambda \in \Omega_{\varepsilon,
C_{\varepsilon}}\}. $

4. As in \cite{k92}, we define $\Omega^{(M_p)}$ as a subset of
${\mathbb C}$ which contains a domain of the form
$$\Omega^{M_p}_{k,C} :=\{\lambda \in {\mathbb C} : Re \lambda \geq
M(k|\lambda|)+ C\},$$
 for some  $k>0$  and  $C>0,$ in
the Beurling case,

\noindent resp., $\Omega^{\{M_p\}}$ as a subset of ${\mathbb C}$
which contains a domain of the form
$$\Omega^{M_p}_{k,C_k} :=\{\lambda \in {\mathbb C} : Re \lambda \geq
M(k|\lambda|)+ C_k\},$$ for every $k>0$ and  the corresponding
$C_k>0,$ in the Roumieu case.

\noindent We use the notation $\Omega^*$ for the common case and put
$(\Omega^{*})^2:=\{\lambda^2 : \lambda  \in \Omega^*\}$.
We define $(\Omega^{M_p}_{k,C})^2$ and $(\Omega^{M_p}_{k,C_k})^2$ in an adequate way. 

5. As in \cite{me152} (cf. also J. Chazarain \cite{cha}), we use the
ultra-logarithmic regions 
%
$$
\Lambda_{\alpha,\beta,\gamma}:=\{\lambda \in {\mathbb C}: Re\lambda
\geq \frac{M(\alpha \lambda)}{\gamma}+\beta\},\quad \alpha,\ \beta,\
\gamma
>0
$$
and  define $\Lambda^{2}_{\alpha,\beta,\gamma}:=\{\lambda^{2}:
\lambda \in \Lambda_{\alpha,\beta,\gamma}\}$. Note that $(M.2)$
implies that, for every $\alpha,\ \beta,\ \gamma
>0,$ there exist  $\alpha'>0$ and
$\beta'>0$ so that $\Lambda_{\alpha',\beta',1} \subset
\Lambda_{\alpha,\beta,\gamma}.$ \vspace{0.1cm}

6. Let $\alpha,\ \beta>0$. The exponential region $E(\alpha,\beta)$
is defined in \cite{a22} by
$$
E(\alpha, \beta):=\{\lambda:\; Re\lambda\geq \beta,\mbox{
}|Im\lambda|\leq e^{\alpha Re\lambda}\};\ E^{2}(\alpha,
\beta):=\{\lambda^{2}: \lambda \in E(\alpha, \beta)\}.
$$ 

7. Let $0<\alpha \leq \pi.$ Then
 $\Sigma_{\alpha}:=\{re^{i \theta}:\mbox{   }r>0,\mbox{
}|\theta|<\alpha\}. $

8. We use occasionally the following condition for $K$:
\begin{itemize}
\item[(P1)]  $K\in L^1_{loc}([0,\infty))$
 is Laplace transformable, i.e.,  there exists $\beta \in {\mathbb R}$ so that
$\tilde{K}(\lambda)={\mathcal L}(K)(\lambda):=\int
\limits^{\infty}_{0}e^{-\lambda t}K(t)dt$ exists for all $\lambda
\in {\mathbb C}$ with $Re\lambda>\beta.$
\end{itemize}

\noindent Put abs$(K):=$inf$\{
Re\lambda : \tilde{K}(\lambda) \mbox{ exists} \}.$ 
In some statements, the next condition is required for $K$
satisfying (P1): 
\begin{itemize}
\item[(P2)] $\tilde{K}(\lambda ) \neq 0,\; Re\lambda >\beta $, where
$\beta \geq \mbox{ abs}(K).$
\end{itemize}
In general, (P2) does not hold for  exponentially bounded functions,
cf. ~\cite[Theorem 1.11.1]{a43}.

9. A function $K\in L^1_{loc}([0,\tau)),  \;\tau>0,$ is called a
kernel if for every $\phi \in C([0,\tau))$, the assumption $\int
\limits^{t}_{0}K(t-s)\phi(s)ds=0,\ t\in [0,\tau)$, implies $\phi
\equiv 0.$
According to Titchmarsh's theorem, $K$ is a kernel if $0\in$supp$K.$

For the later use we recall a family of kernels, see ~\cite[p.
107]{a43}:
$$K_{\delta}(t):=\frac{1}{2\pi i}\int \limits^{r+i\infty}_{r-i\infty}
e^{\lambda t-\lambda^{\delta}}d\lambda,\; t\geq 0,\; 0<\delta<1,\;
r>0, \mbox{ where } 1^{\delta}=1.$$ Note,
$K_{1/2}(t)=\frac{1}{2\sqrt{\pi t^{3}}}e^{-\frac{1}{4t}},\ t>0
\;(K_{1/2}(0)=0).$

\section{ $K$-convoluted $C$-cosine functions}

We assume in the sequel that $K$ is not identical to zero. The
definitions of (local) $K$-convoluted $C$-semigroups and
exponentially bounded, $K$-convoluted $C$-semigroups are recalled in
the Appendix.


\begin{defn}\label{1.17}
Let $A$ be a closed operator and $K \in L^{1}_{loc}([0,\tau))$, $
0<\tau \leq \infty $. If there exists a strongly continuous operator family $%
(C_{K}(t))_{t\in \lbrack 0,\tau )}$ such that:

\begin{itemize}
\item[(i)] $C_{K}(t)A\subset AC_{K}(t)$, $t\in [0,\tau )$,

\item[(ii)] $C_{K}(t)C=CC_{K}(t),\ t\in [0,\tau)\mbox{ and}$

\item[(iii)] $
\int\limits_{0}^{t}(t-s)C_{K}(s)xds\in D(A),\; x\in E,\; t\in
[0,\tau ) \;\mbox{ and }$
\begin{equation}\label{cicika}
A\int\limits_{0}^{t}(t-s)C_{K}(s)xds=C_{K}(t)x-\Theta (t)Cx,\mbox{
where } \Theta (t) := \int \limits^{t}_{0} K(s)ds,\;
\end{equation}
\end{itemize}
then it is said that $A$ is a subgenerator of a $K$-convoluted
$C$-cosine function $(C_{K}(t))_{t\in [0,\tau)}.$ If $\tau=\infty$,
then we say that $(C_{K}(t))_{t\geq 0}$ is an exponentially bounded,
$K$-convoluted $C$-cosine function with a subgenerator $A$ if,
additionally, there exist $M>0$ and $\omega \in {\mathbb R}$ such
that $||C_{K}(t)||\leq Me^{\omega t},\ t\geq 0.$
\end{defn}

As a consequence of (i) and (iii), we have $CA\subset AC.$ Indeed,
if $x\in D(A),$ choose a $t\in [0,\tau)$ with $\Theta(t)\neq 0.$
Then (i) and (iii) implies $C_{K}(t)Ax-\Theta(t)CAx=A\int
\limits^{t}_{0}(t-s)C_{K}(s)Axds=A^{2}\int
\limits^{t}_{0}(t-s)C_{K}(s)xds=A[C_{K}(t)x-\Theta(t)Cx].$ Since
$C_{K}(t)x\in D(A),$ we obtain $Cx\in D(A)$ and $CAx=ACx.$

Put in Definition \ref{1.17},
$K(t) = \frac{t^{\alpha-1}}{\Gamma (\alpha)}, \ t\in [0,\tau),\;
\alpha >0.$ Then
$(C_{K}(t))_{ t\in [0,\tau )} $ is an
$\alpha$-times integrated $C$-cosine function. We point out that C.
Lizama used in \cite{liz121} a slight modification of (\ref{cicika})
and (\ref{ci}) (see Appendix) in the case of $\alpha$-times
integrated cosine functions and semigroups.

The integral generator of $(C_{K}(t))_{t\in [0,\tau)}$ is defined by
\\ $\{(x,y)\in
E^{2} : C_{K}(t)x-\Theta(t)Cx=\int \limits^{t}_{0}(t-s)C_{K}(s)yds,\
t\in[0,\tau)\}.$  The integral generator of  $(C_{K}(t))_{t\in
[0,\tau)}$ is a closed linear operator which is an extension of any
subgenerator of $(C_{K}(t))_{t\in [0,\tau)}$. Even if $(C(t))_{t\geq
0}$ is a global, exponentially bounded $C$-cosine function, the set
of all subgenerators of $(C(t))_{t\geq 0}$ need not be monomial.
This can be viewed by transferring  ~\cite[Example 2.13]{w999} to
cosine functions. Moreover, the set of all subgenerators of a
$K$-convoluted $C$-cosine function can have infinitely many
elements. In order to illustrate this fact, choose an arbitrary
$K\in L^{1}_{loc}([0,\infty)).$ Put $E:=l_{\infty},$ $C\langle x_{n}
\rangle :=\langle 0,x_{1},0,x_{2},0,x_{3}, \ldots \rangle$ and
$C_{K}(t)\langle x_{n} \rangle:=\Theta(t)C\langle x_{n} \rangle,\
t\geq 0,\ \langle x_{n} \rangle \in E.$ If $I \subset 2{\mathbb
N}+1,$ define $E_{I}:=\{ \langle x_{n} \rangle \in E : x_{i}=0,
\mbox{ for all} \ i\in (2{\mathbb N}+1) \setminus I \}.$ Then
$E_{I}$ is a closed subspace of $E$ which contains $R(C).$ Clearly,
$E_{I_{1}}\neq E_{I_{2}},$ if $I_{1}\neq I_{2}.$ Define a closed
linear operator $A_{I}$ on $E$ by: $D(A_{I})=E_{I}$ and
$A_{I}\langle x_{n} \rangle=0, \ \langle x_{n} \rangle \in
D(A_{I}).$ It is straightforward to see that every subgenerator of
$(C_{K}(t))_{t\geq 0}$ is of the form $A_{I},$ for some $I \subset
2{\mathbb N}+1.$ Hence, in this example, there exist the continuum
many subgenerators of $(C_{K}(t))_{t\geq 0}.$ See also
~\cite[Example 2.14]{w999} for a more complicated construction in
the case of global $C$-semigroups.

 If $C=I,$ then  the proof of ~\cite[Proposition 2.2]{mn}, with
slight modifications, shows that every subgenerator of a (local)
$K$-convoluted cosine function $(C_{K}(t))_{t\in [0,\tau)}$
coincides with the integral generator of $(C_{K}(t))_{t\in
[0,\tau)}.$

{\bf Open problem.} The authors do not know whether the set of all
subgenerators $(C_{K}(t))_{t\in [0,\tau)}$ must be monomial if
$C\neq I$ and $\overline{R(C)}=E.$

We need the following useful extension of ~\cite[Proposition
1.3]{ko99}.

\begin{prop}\label{pro}
Let $A$ be a closed operator and let $K\in L^{1}_{loc}([0,\tau))$,
$0<\tau \leq \infty $. Then the following assertions are equivalent:

\begin{itemize}
\item[(a)] $A$ is a subgenerator of a $K$-convoluted $C$-cosine function $(C_{K}(t))_{t\in [0,\tau )}$ in $E$.

\item[ (b)] The operator ${\mathcal A}\equiv $ $(
\begin{array}{cc}
0 & I \\
A & 0%
\end{array}%
) $ \ is a subgenerator of a $\Theta$-convoluted ${\mathcal
C}$-semigroup $(S_{\Theta}(t))_{t\in \lbrack 0,\tau )}$ in $E^{2},$
where ${\mathcal C}:=(\begin{array}{cc}
C & 0 \\
0 & C
\end{array}
).$ In this case:
\end{itemize}

\[
S_{\Theta}(t)=\left(
\begin{array}{cc}
\int\limits_{0}^{t}C_{K}(s)ds & \int\limits_{0}^{t}(t-s)C_{K
}(s)ds \\
C_{K}(t)-\Theta(t)C\; &
\int\limits_{0}^{t}C_{K}(s)ds%
\end{array}%
\right) ,\;\; 0\leq t<\tau,
\]
and the integral generators of $(C_{K}(t))_{t\in [0,\tau )}$ and
$(S_{\Theta}(t))_{t\in [ 0,\tau )},$ denoted respectively by $B$ and
${\mathcal B},$ satisfy ${\mathcal B}=(\begin{array}{cc}
0 & I \\
B & 0%
\end{array}).$ Furthermore,
if $K$ is a kernel, then the integral generator of $(C_{K}(t))_{t\in
[0,\tau )},$ resp., $(S_{\Theta}(t))_{t\in \lbrack 0,\tau )}$ is
$C^{-1}AC,$ resp., ${\mathcal C}^{-1}{\mathcal A}{\mathcal C} \equiv
\left(
\begin{array}{cc}
0 & I \\
C^{-1}AC\; & 0
\end{array} \right). $
\end{prop}

{\em Proof:}
 (a) $\Rightarrow $ (b) The properties of $(C_{K}(t))_{t\in
[0,\tau )}$ and $CA\subset AC$ imply
that $%
(S_{\Theta}(t))_{t\in \lbrack 0,\tau )}$ is a strongly continuous
operator family in $E^{2}$ which satisfies $S_{\Theta}(t){\mathcal A}\subset {\mathcal %
A}S_{\Theta}(t)$ and $S_{\Theta}(t){\mathcal C}= {\mathcal
C}S_{\Theta }(t)$ for $0\leq t<\tau .$ Furthermore,
\[
{\mathcal A}\int\limits_{0}^{t}S_{\Theta}(s){x\choose y}ds=
{\mathcal A}%
\int\limits_{0}^{t}\left(
\begin{array}{cc}
\int\limits_{0}^{s}C_{K}(r)xdr+\int\limits_{0}^{s}(s-r)C_{K}(r)ydr \\
C_{K}(s)x-\Theta (s) Cx+\int\limits_{0}^{s}C_{K}(r)ydr
\end{array}\right) ds.
\]
\[
={\mathcal A}\left(
\begin{array}{cc}
\int\limits_{0}^{t}(t-s)C_{K}(s)xds+\int\limits_{0}^{t}%
\frac{(t-s)^{2}}{2}C_{K}(s)yds \\
\int\limits_{0}^{t}C_{K}(s)xds-\int \limits^{t}_{0}\Theta
(s)Cxds+\int\limits_{0}^{t}(t-s)C_{K}(s)yds
\end{array}
\right)
\]%
$$
=\left( \begin{array}{cc} \int\limits_{0}^{t}C_{K}(s)xds-\int
\limits^{t}_{0}\Theta (s)Cxds+\int\limits_{0}^{t}(t-s)C_{K }(s)yds
\\
C_{K}(t)x-\Theta (t)Cx+\int\limits_{0}^{t}C_{K }(s)yds-\int
\limits^{t}_{0}\Theta (s)Cyds
\end{array}
\right)
$$
$$
=S_{\Theta}(t){x\choose y} -\int \limits^{t}_{0}\Theta
(s){{Cx}\choose{Cy}}ds, \; 0\leq t<\tau .\; 
$$

(b) $\Rightarrow $ (a) Put $S_{\Theta}(t)=\left(
\begin{array}{cc}
S_{\Theta}^{1}(t) & S_{\Theta}^{2}(t) \\
S_{\Theta}^{3}(t) & S_{\Theta}^{4}(t)%
\end{array}%
\right), \; {t\in \lbrack 0,\tau )},$ where \newline
$S_{\Theta}^{i}(t)\in L(E)$, $i\in \{1,2,3,4\}$, $0\leq t<\tau $. A
simple consequence of $S_{\theta}(t){\mathcal C}={\mathcal
C}S_{\Theta}(t),\ t\in [0,\tau)$ is:
$S_{\Theta}^{i}(t)C=CS_{\Theta}^{i}(t),\ t\in [0,\tau),\ i\in
\{1,2,3,4\}.$ Since $ S_{\Theta}{\mathcal A}\subset {\mathcal
A}S_{\Theta}:$
\[
\begin{array}{c}
S_{\Theta}^{1}(t)x+S_{\Theta}^{2}(t)y\in D(A), \\[0.2cm]
S_{\Theta}^{1}(t)y+S_{\Theta}^{2}(t)Ax=S_{\Theta}^{3}(t)x+S_{\Theta}^{4}(t)y, \\[0.2cm]
S_{\Theta}^{3}(t)y+S_{\Theta}^{4}(t)Ax=A(S_{\Theta}^{1}(t)x+S_{\Theta}^{2}(t)y),\;\; 0\leq t<\tau ,\;\; x\in D(A),%
\;\; y\in E.%
\end{array}%
\]%
Hence, $S_{\Theta}^{3}(t)x=S_{\Theta}^{2}(t)Ax$, $x\in D(A)$, and $%
S_{\Theta}^{3}(t)y=AS_{\Theta}^{2}(t)y$, $y\in E$, $0\leq t<\tau $.
This implies that for every $x\in D(A)$, we have
$S_{\Theta}^{3}(t)Ax=AS_{\Theta}^{2}(t)Ax=AS_{\Theta}^{3}(t)x,\; t\in [0,\tau)$. Thus, $%
S_{\Theta}^{3}(t)A\subset AS_{\Theta}^{3}(t)$, $t\in \left[ 0,\tau
\right) $, and $(S_{\Theta}^{3}(t)+\Theta (t)C)_{t\in \lbrack 0,\tau
)}$ is a strongly continuous operator family in $E$. Now, the simple
calculation deduced from
$
{\mathcal A}\int\limits_{0}^{t}S_{\Theta}(s) {x\choose y} ds=S_{\Theta}(t)%
{x\choose y}-\int \limits^{t}_{0}\Theta (s){{Cx}\choose {Cy}}ds, $
gives
\[
\int\limits_{0}^{t}S_{\Theta}^{3}(s)xds+\int\limits_{0}^{t}S_{\Theta}^{4}(s)yds=S_{\Theta}^{1}(t)x+S_{\Theta}^{2}(t)y-\int
\limits^{t}_{0}\Theta (s)C xds \ \mbox{ and}
\]%
\[
A\left[
\int\limits_{0}^{t}S_{\Theta}^{1}(s)xds+\int\limits_{0}^{t}S_{\Theta}^{2}(s)yds\right]
=S_{\Theta}^{3}(t)x+S_{\Theta}^{4}(t)y-\int \limits^{t}_{0}\Theta
(s)Cyds,
\]
for all $0\leq t<\tau ,\;\; x,\ y\in E$. Hence, $\int\limits_{0}^{t}S_{\Theta}^{3}(s)xds=S_{\Theta}^{1}(t)x-%
\int \limits^{t}_{0}\Theta (s)Cxds$ and $A\int\limits_{0}^{t}S_{%
\Theta}^{1}(s)xds=S_{\Theta}^{3}(t)x$, $0\leq t<\tau $, $x\in E$.
Consequently,
\[
A\left[ \int\limits_{0}^{t}(t-s)(S_{\Theta}^{3}(s)x+\Theta
(s)Cx)ds\right]
=A\left[ \int\limits_{0}^{t}(t-s)(\frac{d}{dv}%
S_{\Theta}^{1}(v)x)_{v=s}ds\right]
\]
\[
=A\int\limits_{0}^{t}S_{\Theta}^{1}(s)xds =\left[
S_{\Theta}^{3}(t)x+\Theta (t)Cx\right] -\Theta (t)Cx\;\; ,\mbox{
}0\leq t<\tau ,\;\; x\in E.
\]%
Thus, we have proved that $A$ is a subgenerator of the $K$-convoluted $C$-cosine function $%
(S_{\Theta}^{3}(t)+\Theta (t)C)_{t\in \lbrack
0,\tau )}$. Clearly, $S_{\Theta}^{1}(t)=S_{\Theta}^{4}(t)$ and $%
S_{\Theta}^{2}(t)=\int\limits_{0}^{t}S_{\Theta}^{1}(s)ds$, $0\leq
t<\tau $. 
Next, we will prove that ${\mathcal B}=(\begin{array}{cc}
0 & I \\
B & 0%
\end{array}).$ To see this, fix some $x,\ y,\ x_{1},\ y_{1} \in E.$ Then
$$
S_{\Theta}(t){{x}\choose {y}}-\int
\limits^{t}_{0}\Theta(s){{Cx}\choose {Cy}}ds=\int
\limits^{t}_{0}S_{\Theta}(s){{x_{1}}\choose {y_{1}}}ds,\; t\in
[0,\tau),
$$
iff
$$
C_{K}(t)x-\Theta(t)Cx=\int
\limits^{t}_{0}(t-s)C_{K}(s)y_{1}ds,\mbox{ for all }t\in [0,\tau),
\mbox{ and }y=x_{1}.
$$
Namely, if $S_{\Theta}(t){{x}\choose {y}}-\int
\limits^{t}_{0}\Theta(s){{Cx}\choose {Cy}}ds=\int
\limits^{t}_{0}S_{\Theta}(s){{x_{1}}\choose {y_{1}}}ds,$ $t\in
[0,\tau),$ then
$$
\int \limits^{t}_{0}C_{K}(s)xds+\int
\limits^{t}_{0}(t-s)C_{K}(s)yds-\int \limits^{t}_{0}\Theta(s)Cxds$$
\begin{equation}\label{pirlo}
=\int \limits^{t}_{0}(t-s)C_{K}(s)x_{1}ds+\int
\limits^{t}_{0}\frac{(t-s)^{2}}{2}C_{K}(s)y_{1}ds,\mbox{  and  }
\end{equation}
$$
C_{K}(t)x-\Theta(t)Cx+\int \limits^{t}_{0}C_{K}(s)yds-\int
\limits^{t}_{0}\Theta(s)Cyds$$
\begin{equation}\label{pirlo1}
=\int \limits^{t}_{0}C_{K}(s)x_{1}ds-\int
\limits^{t}_{0}\Theta(s)Cx_{1}ds+\int
\limits^{t}_{0}(t-s)C_{K}(s)y_{1}ds.
\end{equation}
Differentiating (\ref{pirlo}) with respect to $t$, one obtains
$$
C_{K}(t)x+\int \limits^{t}_{0}C_{K}(s)yds-\Theta(t)Cx=\int
\limits^{t}_{0}C_{K}(s)x_{1}ds+\int
\limits^{t}_{0}(t-s)C_{K}(s)y_{1}ds.
$$
The last equality and (\ref{pirlo1}) imply $\int
\limits^{t}_{0}\Theta(s)Cyds=\int
\limits^{t}_{0}\Theta(s)Cx_{1}ds$; consequently, $y=x_{1}.$ Then
(\ref{pirlo1}) gives $ C_{K}(t)x-\Theta(t)Cx=\int
\limits^{t}_{0}(t-s)C_{K}(s)y_{1}ds, t\in [0,\tau)$ and
$(x,y_{1})\in B.$ Conversely,  suppose  that $y=x_{1}$ and that
$(x,y_{1})\in B.$ Then \\ $ C_{K}(t)x-\Theta(t)Cx=\int
\limits^{t}_{0}(t-s)C_{K}(s)y_{1}ds,\; t\in [0,\tau).$ This
implies (\ref{pirlo1}). Integrating (\ref{pirlo1}) with respect to
$t$ one obtains (\ref{pirlo}) and this gives ${{x}\choose
{y}},{{x_{1}}\choose {y_{1}}}\in {\mathcal B}.$ Hence, we have
proved ${{x}\choose {y}},{{x_{1}}\choose {y_{1}}}\in {\mathcal B}$
iff $y=x_{1}$ and $(x, y_{1})\in B.$ Further on, our assumption
$CA\subset AC$ implies ${\mathcal C}{\mathcal A}\subset {\mathcal
A}{\mathcal C}$ and one can employ Proposition \ref{natalija}
given below in order to see that the integral generator of
$(S_{\Theta}(t))_{t\in [0,\tau)}$ is ${\mathcal C}^{-1}{\mathcal
A}{\mathcal C}.$ As a matter of routine, we obtain ${\mathcal
C}^{-1}{\mathcal A}{\mathcal C}=\left(
\begin{array}{cc}
0 & I \\
C^{-1}AC\; & 0
\end{array} \right).$ By the previously given arguments, we know that this implies that the integral generator
of $(C_{K}(t))_{t\in[0,\tau)}$ is $C^{-1}AC.$

\begin{rem}\label{zvonecet}
When $\tau =\infty $ and $\Theta$ is an exponentially bounded
function, then $(C_{K}(t))_{t\geq 0}$ is exponentially bounded if
and only if $(S_{\Theta}(t))_{t\geq 0}$ is exponentially bounded.
\end{rem}

Proposition \ref{pro} implies the following facts which remain true
in the case of convoluted $C$-semigroups (see Appendix).

 Suppose in
this paragraph that $0<\tau \leq \infty$ and that $K\in
L_{loc}^{1}([0,\tau))$ is a kernel. If $A$ and $B$ are
subgenerators of $%
(C_{K}(t))_{t\in \lbrack 0,\tau )},$ then: $C^{-1}AC=C^{-1}BC,$
$C(D(A))\subset D(B)$ and $A=B \Leftrightarrow D(A)=D(B).$ Moreover,
if $A$ is the integral generator of $%
(C_{K}(t))_{t\in \lbrack 0,\tau )},$ then
it can be easily seen that the set of all subgenerators of $%
(C_{K}(t))_{t\in \lbrack 0,\tau )}$ is singleton if $C(D(A))$ is a
core for $D(A),$ cf. also ~\cite[Proposition 2.8]{w999} and
\cite{l1}. It can be proved that all subgenerators of
$(C_{K}(t))_{t\in \lbrack 0,\tau )}$ form a lattice. Further
analysis of such a lattice can be found  \cite{w999}.

If $K$ is kernel, then by Proposition \ref{pro} and the
corresponding statement in the case of semigroups that every (local)
$K$-convoluted $C$-cosine function is uniquely determined by one of
its subgenerators. The standard proof is omitted.

%

Let $(C_{K}(t))_{t\in [0,\tau)}$ be a (local) $K$-convoluted
$C$-cosine function whose integral generator is $A.$ Proposition
\ref{pro} and Proposition \ref{natalija} of Appendix immediately
yield $A=C^{-1}AC.$

Let $\lambda \in {\mathbb C}$ and let $E^{2}$ be endowed  by the
graph norm $||(x,y)||=||x||+||y||,\ x,\ y\in E.$ Then $\lambda^{2}
\in \rho(A)$ iff $\lambda \in \rho({\mathcal A}),$ and, in this
case, we have:
$$||R(\lambda^{2} : A)||\leq ||R(\lambda : {\mathcal
A}) ||,$$
$$||R(\lambda : {\mathcal A})||\leq
(1+|\lambda|)\sqrt{1+|\lambda|^{2}}||R(\lambda^{2} : A)||+1,$$
$$
R(\lambda:{\mathcal A}){{x}\choose {y}}={{R(\lambda^{2}:A)(\lambda
x+y)}\choose {AR(\lambda^{2}:A)x+\lambda
R(\lambda^{2}:A)y}},\mbox{ }x,\mbox{ }y\in E,
$$
see ~\cite[Lemma 1.10]{ko99} for the proof. Further on, $D({\mathcal A}^{n})=D(A^{\lceil \frac{%
n}{2}\rceil })\times D(A^{\lfloor \frac{n}{2}\rfloor })$, $n\in {{\mathbb N}_{0}}$%
\ and $D_{\infty }({\mathcal A})=D_{\infty }(A)\times D_{\infty
}(A)$. Here, $\lfloor t \rfloor:=\sup \{k \in \mathbb{Z} : k\leq
t\}$ and $\lceil t \rceil:=\inf\{k \in \mathbb{Z} : k\geq t\},\ t\in
{\mathbb R}.$

Let us recall (\cite{ku101}) that a closed linear operator $A$ is
stationary dense if $ n(A):=\inf \{ k\in {\mathbb N}_{0} : (\forall
n\geq k)\;\; D(A^{n})\subset \;\;\overline{D{\Bbb (}A^{n+1}{\Bbb
)}}\}<\infty {\Bbb .}$ We will prove that every generator of an
integrated cosine function is stationary dense. In Example
\ref{Beals} we will show that this is not automatically satisfied if
$A$ generates a convoluted cosine function.

\begin{lem}\label{simpli}
Let $A$ be a closed operator. Then $A$ is stationary dense if and
only if ${\mathcal A}$ is stationary dense. Moreover, $n({\mathcal
A}) = 2n(A)$.
\end{lem}
{\em Proof:}
 Assume that $A$ is stationary dense and that $n(A)=n\in
{\mathbb{N}}_{0}$. Let us prove that $D({\mathcal A}^{m})\subset
\overline{D({\mathcal A}^{m+1})},$ for all $m\in {\mathbb{N}}_{0}$
with $m\geq 2n$. Suppose $m=2i$, for some $i\geq n$. We have to
prove that $D(A^{i}) \times D(A^{i}) \subset \overline{D(A^{i+1})
\times D(A^{i})}$. But, this is a consequence of $D(A^{i})\subset
\overline{D(A^{i+1})}$. If $m=2i+1$ for some $i\geq n$, then
$D({\mathcal A}^{m})\subset \overline{D({\mathcal A}^{m+1})}$ is
equivalent with $D(A^{i+1}) \times D(A^{i}) \subset
\overline{D(A^{i+1}) \times D(A^{i+1})}$, which is valid since
$i\geq n$. Thus, ${\mathcal A}$ is stationary dense and $n({\mathcal
A})\leq 2 n(A).$ Furthermore, $n({\mathcal A})=0,$ if $n(A)=0.$
Suppose $n({\mathcal A})< 2 n(A).$ If $ n({\mathcal A})=2i$, for
some $i\in \{ 0, 1, \ldots, n-1 \}$, then $D(A^{i}) \times D(A^{i})
\subset \overline{D(A^{i+1}) \times D(A^{i})}.$ It gives
$D(A^{i})\subset \overline{D(A^{i+1})}$ and the contradiction is
obvious. Similarly, if $ n({\mathcal A})=2i+1$, for some $i\in \{ 0,
1, \ldots, n-1 \}$, then $D(A^{i+1}) \times D(A^{i}) \subset
\overline{D(A^{i+1}) \times D(A^{i+1})}.$ Again, $D(A^{i})\subset
\overline{D(A^{i+1})}$ and this is in contradiction with $n(A)=n$.
Hence, we have proved that ${\mathcal A}$ is stationary dense and $
n({\mathcal A})= 2 n(A).$ Assume conversely that ${\mathcal A}$ is
stationary dense. Similarly as in the first part of the proof, one
obtains that $A$ is stationary dense. Then we know that $n({\mathcal
A})=2n(A).$

\begin{prop}\label{stacija}
Let $A$ be the generator of an $\alpha$-times integrated cosine
function $(C_{\alpha}(t))_{t\in [0,\tau )}$, $0<\tau \leq \infty $,
$\alpha >0.$ Then $n(A)\leq \lfloor \frac{\lceil \alpha
\rceil+1}{2}\rfloor $.
\end{prop}
{\em Proof:}
Due to Proposition \ref{pro} and ~\cite[Proposition 2.4(a)]{l114},
the operator ${\mathcal A}$ is the generator of an $%
(\lceil \alpha \rceil+1)$-times integrated semigroup $(S_{\lceil
\alpha \rceil+1}(t))_{t\in \lbrack 0,\tau )}$. Thus, an application
of ~\cite[Corollary 1.8]{ku101} gives $n({\mathcal A})\leq \lceil
\alpha \rceil+1 $. Now the proposition follows from Lemma
\ref{simpli}.

{\sc Comment and Problem.} As it is illustrated in ~\cite[Example
3.15.5, p. 224]{a43}, the generator $B$ of the standard translation
group on $L^{1}({\mathbb R})$ fulfills the next statement:
$A:=(B^{\ast})^{2}$ is the non-densely defined generator of a sine
function in $L^{\infty}({\mathbb R}).$ Proposition \ref{stacija}
implies $n(A)=1.$ Hence, in the general situation of previous
proposition, the estimate $n(A)\leq \lfloor \frac{\lceil \alpha
\rceil +\beta}{2}\rfloor $ cannot be proved for any $\beta \in
[0,1)$ since here $n(A)=1$ and $\alpha=1.$ The next problem can be
posed: Given an arbitrary $\alpha>0,$ is it possible to construct a
Banach space $E_{\alpha},$ a closed linear operator $A_{\alpha}$ on
$E_{\alpha}$ which generates a (local) $\alpha$-times integrated
cosine function and satisfies $n(A_{\alpha})=\lfloor \frac{\lceil
\alpha \rceil+1}{2}\rfloor ?$


\section{Global exponentially bounded $K$-convoluted $C$-cosine functions}

Recall, the $C$-resolvent set of $A$, denoted by $\rho_{C}(A),$ is
defined by \\ $\rho_{C}(A):=\{ \lambda \in {\mathbb C} : R(C)\subset
R(\lambda-A)\mbox{ and } \lambda-A \mbox{ is injective} \}.$

\begin{thm}\label{bau}
Assume that $K$ satisfies (P1) and that $A$ is a closed linear
operator.
\begin{itemize}
\item[(a)] Assume that $A$ is a subgenerator of an exponentially
bounded, $K$-convoluted $C$-cosine function $(C_{K}(t))_{t\geq 0}$
and $||C_{K}(t)||\leq Me^{\omega t},\ t\geq 0$, for some $M>0$ and
$\omega \geq 0.$ If $\omega_{1}=\max(\omega, \beta)$, then
\begin{equation}\label{bekaa}
\{\lambda^{2} : \ Re \lambda
>\omega_{1},\mbox{   }\tilde{K}(\lambda)\neq 0 \} \subset \rho_{C}
(A)\mbox{ and}
\end{equation}
\begin{equation}\label{beko}
\lambda (\lambda ^{2}-A)^{-1}Cx=\frac{1}{\tilde{K}(\lambda )}%
\int\limits_{0}^{\infty }e^{-\lambda t}C_{K}(t)xdt,\; x\in E, \; Re
\lambda > \omega_{1},\mbox{ }\tilde{K}(\lambda)\neq 0.\;\;
\end{equation}
\item[(b)] Suppose that $(C_{K}(t))_{t\geq 0}$ is a strongly continuous
operator family satisfying $||C_{K}(t)||\leq Me^{\omega t},$ $t\geq
0,\ \omega \geq 0.$ Put $\omega_{1}=\max(\omega, \beta).$ If
(\ref{bekaa}) and (\ref{beko}) are fulfilled, then
$(C_{K}(t))_{t\geq 0}$ is an exponentially bounded, $K$-convoluted
$C$-cosine function with a subgenerator $A.$
\end{itemize}
\end{thm}
{\em Proof:}
 (a) Fix a $\lambda \in {\mathbb C}$ satisfying
$Re\lambda>\omega_{1}$ and $\tilde{K}(\lambda)\neq 0.$ Since
(\ref{cicika}) is assumed and $A$ is closed, we have
\[
{\mathcal L}(C_{K}(t)x)(\lambda )=\frac{\tilde{K}(\lambda )}{\lambda }Cx+\frac{1%
}{\lambda ^{2}}A{\mathcal L}(C_{K}(t)x)(\lambda ),\; \mbox{i.e.,}
\]%
\begin{equation}\label{cole}
(\lambda ^{2}-A){\mathcal L}(C_{K}(t)x)(\lambda )=\lambda
\tilde{K}(\lambda )Cx,\ x\in E.
\end{equation}
Hence, $R(C)\subset R(\lambda ^{2}-A).$ Let us show that $%
\lambda ^{2}-A$ is injective. Suppose $%
(\lambda ^{2}-A)x=0$. This implies
\[
C_{K}(t)x-\Theta (t)Cx=\int\limits_{0}^{t}(t-s)C_{K}(s)Axds=\lambda
^{2}\int\limits_{0}^{t}(t-s)C_{K}(s)xds,\;t\geq 0,
\]%
and consequently, $$
{\mathcal L}(C_{K}(t)x)(\lambda )=\frac{\tilde{K}(\lambda )}{%
\lambda }Cx+{\mathcal L}(C_{K}(t)x)(\lambda ).$$ Thus, $Cx=0$ and
$x=0$. Then (\ref{cole}) implies $(w_{1}^{2},\infty )\subset \rho
_{C}(A)$ and (\ref{beko}).

(b) Fix an $x\in E$ and a $ \lambda \in {\mathbb C}$ with
$\tilde{K}(\lambda)\neq 0$ and $Re\lambda>\omega_{1}.$ Then
(\ref{beko}) and $CA\subset AC$ imply $(\lambda
^{2}-A)^{-1}C^{2}x=C(\lambda ^{2}-A)^{-1}Cx.$ The previous equality
gives
\[
\frac{1}{\tilde{K}(\lambda )}\int\limits_{0}^{\infty }e^{-\lambda
t}C_{K}(t)Cxdt=\frac{1}{\tilde{K}(\lambda )}\int\limits_{0}^{\infty
}e^{-\lambda t}CC_{K}(t)xdt.
\]
Since $K\neq 0$ in $L^{1}_{\mbox{loc}}([0,\infty)),$ it follows
$$
\{z\in {\mathbb C}: Rez \geq \omega_{1}\}=\overline{\{z\in {\mathbb
C}: Rez>\omega_{1},\mbox{   }\tilde{K}(z)\neq 0\}}.$$ Hence,
$$\int\limits_{0}^{\infty }e^{-\lambda
t}C_{K}(t)Cxdt=\int\limits_{0}^{\infty }e^{-\lambda t}CC_{K}(t)xdt$$
if $Re\lambda>\omega_{1}.$ The uniqueness theorem for the Laplace
transform implies $CC_{K}(t)=C_{K}(t)C,\ t\geq 0.$ Assume now $x\in
D(A)$. If $Re\lambda>\omega_{1}$ and $\tilde{K}(\lambda)\neq 0,$
then
$$\lambda (\lambda ^{2}-A)^{-1}CAx=\frac{1}{\tilde{K}(\lambda )}%
\int\limits_{0}^{\infty }e^{-\lambda t}C_{K}(t)Axdt,\;\mbox{i.e.,}
$$
$$\lambda A(\lambda ^{2}-A)^{-1}Cx=\frac{1}{\tilde{K}(\lambda )}%
\int\limits_{0}^{\infty }e^{-\lambda t}C_{K}(t)Axdt.
$$
A consequence is $A\int\limits_{0}^{\infty }e^{-\lambda
t}C_{K}(t)xdt=\int\limits_{0}^{\infty }e^{-\lambda t}C_{K}(t)Axdt.$
Using the closedness of $A$ and the same arguments as above, we
obtain that the last equality holds for every $\lambda \in {\mathbb
C}$ with $Re\lambda>\omega_{1}.$ Now one can apply
~\cite[Proposition 1.7.6]{a43} in order to conclude that
$C_{K}(t)A\subset AC_{K}(t),\ t\geq 0.$ Let $Re\lambda>\omega_{1}$
and $\tilde{K}(\lambda)\neq 0.$ We have
$${\mathcal L}\left(
\int\limits_{0}^{t}(t-s)C_{K}(s)xds\right) (\lambda ) = {\mathcal
L}(t)(\lambda ){\mathcal L}(C_{K}(t)x)(\lambda )$$
$$ =\frac{1}{\lambda
^{2}}\tilde{K}(\lambda )\lambda (\lambda ^{2}-A)^{-1}Cx
=\frac{\tilde{K}(\lambda )}{\lambda }(\lambda ^{2}-A)^{-1}Cx,$$
which implies
$$A\left( {\mathcal L}\left(
\int\limits_{0}^{t}(t-s)C_{K}(s)xds\right) (\lambda )\right)
=\tilde{K}(\lambda )\lambda (\lambda ^{2}-A)^{-1}Cx-\frac{\tilde{K
}(\lambda )}{\lambda }Cx$$ $$={\mathcal L}(C_{K}(t)x-\Theta
(t)Cx)(\lambda ).$$ Similarly as above, the closedness of $A$ and
the use of ~\cite[Proposition 1.7.6]{a43} imply that (\ref{cicika})
holds.

Note that there exist examples of local integrated $C$-cosine
functions and semigroups whose integral generators have the empty
$C$-resolvent sets (cf. \cite{l114}).

In the next statement, we relate exponentially bounded, convoluted
$C$-semi-groups to exponentially bounded, convoluted $C$-cosine
functions.

\begin{prop}\label{veza}
Let $K$ satisfy (P1). Suppose that $A$ and $-A$ are subgenerators of
exponentially bounded, $K$-convoluted $C$-semigroups. Then $A^{2}$
is a subgenerator of an exponentially bounded, $K$-convoluted
$C$-cosine function.
\end{prop}
{\em Proof:}
Suppose that $A$ and $-A$ are subgenerators of exponentially
bounded, $K$-convoluted $C$-semigroups $(S_{K}(t))_{t\geq 0}$ and
$(V_{K}(t))_{t\geq 0},$ respectively. Define\\
$C_{K}(t):=\frac{1}{2}(S_{K}(t)+V_{K}(t)),\ t\geq 0.$ We will prove
that $A^{2}$ is a subgenerator of a $K$-convoluted $C$-cosine
function $(C_{K}(t))_{t\geq 0}.$ Clearly, $(C_{K}(t))_{t\geq 0}$ is
an exponentially bounded operator family. Arguing similarly as in
the proof of Theorem \ref{bau}, we obtain that there is an
$\omega_{1}>0$ such that
\begin{equation}\label{zeljez}
\{ \lambda \in {\mathbb C} : Re \lambda>\omega_{1}, \
\tilde{K}(\lambda)\neq 0 \} \subset \rho_{C}(A) \cap \rho_{C}(-A)
\end{equation}
and that
$$
(\lambda-A)^{-1}Cx=\frac{1}{\tilde{K}(\lambda )}
\int\limits_{0}^{\infty }e^{-\lambda t}S_{K}(t)xdt,\; x\in E, \; Re
\lambda > \omega_{1},\mbox{ }\tilde{K}(\lambda)\neq 0.
$$
The previous equality poses the natural analog for $-A$ and
$(V_{K}(t))_{t\geq 0}.$ Fix  $\lambda \in {\mathbb C}$ with $Re
\lambda>\omega_{1}$ and $\tilde{K}(\lambda)\neq 0.$ Due to
(\ref{zeljez}), we have that $\lambda^{2}-A^{2}$ is injective.
Moreover, it is straightforward to see that $R(C)\subset
R(\lambda^{2}-A^{2})$ and that
$$
(\lambda^{2}-A^{2})^{-1}Cx=\frac{1}{2\lambda}[(\lambda-A)^{-1}Cx +
(\lambda +A)^{-1}Cx]
$$
$$
=\frac{1}{\lambda \tilde{K}(\lambda)}\int
\limits^{\infty}_{0}e^{-\lambda
t}[\frac{1}{2}(S_{K}(t)x+V_{K}(t)x)]dt =\frac{1}{\lambda
\tilde{K}(\lambda)}\int \limits^{\infty}_{0}e^{-\lambda
t}C_{K}(t)xdt,\ x\in E.
$$
The proof ends an application of Theorem \ref{bau}.

\section{On the abstract Weierstrass
formula}

We will state the convoluted version of the abstract Weierstrass
formula (Theorem \ref{mance}). First we introduce the class of
analytic $K$-convoluted $C$-semigroups.

\begin{defn}\label{anasem}
Let $0<\alpha\leq \frac{\pi}{2}$ and let $A$ be a subgenerator of a
$K$-convoluted $C$-semigroup $(S_{K}(t))_{t\geq 0}.$ Then we say
that $(S_{K}(t))_{t\geq 0}$ is an analytic $K$-convoluted
$C$-semigroup of angle $\alpha$ having $A$ as a subgenerator, if
there exists an analytic function ${\bf S}_{K}:\Sigma_{\alpha} \to
L(E)$ which satisfies
\begin{itemize}
\item[(i)] ${\bf S}_{K}(t)=S_{K}(t),\ t>0,$ and \vspace{0.1cm}
\item[(ii)] $\lim_{z\rightarrow 0, z\in \Sigma_\gamma}{\bf S} _{K}(z)x=0,$
for all $\gamma\in (0,\alpha)$ and $x\in E.$
\end{itemize}
It is said that $A$ is a subgenerator of an exponentially bounded,
analytic $K$-convoluted $C$-semigroup $(S_{K}(t))_{t\geq 0}$ of
angle $\alpha$, if for every $\gamma\in(0,\alpha),$ there exist
$M_{\gamma}>0$ and $\omega_{\gamma}>0$ such that $||{\bf
S}_{K}(z)||\leq M_{\gamma}e^{\omega_{\gamma}Rez},\ z\in
\Sigma_{\gamma}$.
\end{defn}

We also write $S_{K}$ for ${\bf S}_{K}.$ If $C=I,$ the previous
definition has been recently introduced in \cite{mn}. Although one
can reformulate a great part of facts known for analytic convoluted
semigroups in general case, we focus our attention on the next
result which improves ~\cite[Theorem 6.3]{mn}.

\begin{thm}\label{anali}
Assume $0<\alpha \leq \frac{\pi}{2}$, $K$ satisfies (P1) and $\omega
\geq \max(0, \mbox{abs}(K)).$ Suppose that $A$ is a closed linear
operator with\\ $\{ \lambda \in {\mathbb C} : Re\lambda>\omega,\
\tilde{K}(\lambda)\neq 0\} \subset \rho_{C}(A)$ and that the
function
$$\lambda \mapsto \tilde{K}(\lambda)(\lambda - A)^{-1}C,\;
Re\lambda>\omega,\ \tilde{K}(\lambda)\neq 0,$$ can be analytically
extended to a function
$$\tilde{q}:
\omega+\Sigma_{\frac{\pi}{2}+\alpha}\to L(E)$$ satisfying
$$
\sup_{\lambda \in
\omega+\Sigma_{\frac{\pi}{2}+\gamma}}||(\lambda-\omega)\tilde{q}(\lambda)||<\infty,
\; \gamma \in (0, \alpha)\mbox{  and  } \lim_{\lambda \rightarrow
+\infty}\lambda\tilde{q}(\lambda)=0.
$$
Then $A$ is a subgenerator of an exponentially bounded, analytic
$K$-convoluted
\\ $C$-semigroup of angle $\alpha$.
\end{thm}
{\em Proof:}
The use of ~\cite[Theorem 2.6.1]{a43} implies that there
exists an analytic function $S_{K}:\Sigma_{\alpha}\rightarrow L(E)$
so that $ \sup_{z\in \Sigma_{\gamma}}||e^{-\omega
z}S_{K}(z)||<\infty,\mbox{ for all  }$ $\gamma \in (0,\alpha)$ and
that
$$\tilde{q}(\lambda)=\int \limits^{\infty}_{0}e^{-\lambda
t}S_{K}(t)dt,\mbox{ }Re\lambda>\omega. $$ Put $S_{K}(0):=0,$ fix
$x\in E$ and  $\gamma \in (0,\alpha).$ We will prove that
$$
\lim_{z\rightarrow 0, z\in \Sigma_\gamma}S_{K}(z)x=0.
$$
Note that $f(z):=e^{-\omega z}S_{K}(z)x,\ z\in \Sigma_{\alpha}$, is
analytic and $\sup_{z\in \Sigma_{\gamma}}||f(z)||<\infty.$ By
~\cite[Proposition 2.6.3]{a43}, it is enough to show
$\lim_{t\downarrow 0}S_{K}(t)x=0.$ This is a consequence of the
assumption $\lim_{\lambda \rightarrow
+\infty}\lambda\tilde{q}(\lambda)=0$ and a Tauberian type theorem
 ~\cite[Theorem 2.6.4]{a43}. It follows that $(S_{K}(t))_{t\geq
0}$ is a strongly continuous, exponentially bounded operator family
which satisfies
$$
\tilde{K}(\lambda)(\lambda - A)^{-1}Cx=\int
\limits^{\infty}_{0}e^{-\lambda t}S_{K}(t)xdt, \;\lambda \in
{\mathbb C}, Re\lambda>\omega,  \tilde{K}(\lambda)\neq 0.
$$
Similarly as in the proof of Theorem \ref{bau} (cf. also \cite{koj}
and \cite{mn}), we have that $A$ is a subgenerator of an
exponentially bounded, $K$-convoluted $C$-semigroup
$(S_{K}(t))_{t\geq 0}.$ Since $(S_{K}(t))_{t\geq 0}$ verifies
conditions (i) and (ii), given in  Definition $\ref{anasem},$
$(S_{K}(t))_{t\geq 0}$ is an exponentially bounded analytic
$K$-convoluted $C$-semigroup of angle $\alpha$ having $A$ as a
subgenerator.

The main result of this section reads as follows.

\begin{thm}\label{mance}
Assume that for some $M>0$ and $\beta>0:$ $|K(t)|\leq Me^{\beta t},\
t\geq 0.$ Let $A$ be a subgenerator of an exponentially bounded
$K$-convoluted $C$-cosine function $(C_{K}(t))_{t\geq 0}.$ Then $A$
is a subgenerator of an exponentially bounded analytic
$K_{1}$-convoluted $C$-semigroup $(S(t))_{t\geq 0}$ of angle
$\frac{\pi}{2}$, where:
$$
K_{1}(t):=\int \limits^{\infty}_{0}\frac{se^{-s^{2}/4t}}{2\sqrt{\pi}
t^{3/2}}K(s)ds \mbox{   and   } S(t)x:=\frac{1}{\sqrt{\pi t}}\int
\limits^{\infty}_{0}e^{-s^{2}/4t}C_{K}(s)xds,\ t>0,\ x\in E.
$$
\end{thm}
{\em Proof:}
We follow the proof of the abstract Weierstrass formula (cf.
~\cite[p. 220]{a43}). Due to ~\cite[Proposition 1.6.8]{a43}, $K_{1}$
fulfills (P1), abs$(K_{1})\geq \beta ^{2}$ and
$\widetilde{K_{1}}(\lambda)=\tilde{K}(\sqrt{\lambda}), $ $
Re\lambda>\beta^{2}$. 
Let $x\in E$ be fixed. Putting $r=s/\sqrt{t}$, and using the
dominated convergence theorem after that, one obtains
\begin{equation}\label{gruja}
S(t)x=\int
\limits^{\infty}_{0}\frac{e^{-r^{2}/4}}{\sqrt{\pi}}C_{K}(r\sqrt{t})xdr
\rightarrow 0,\ t\rightarrow 0+.
\end{equation}
Define $S(0):=0$. By (\ref{gruja}), $(S(t))_{t\geq 0}$ is a strongly
continuous, exponentially bounded operator family. Furthermore, one
can employ Theorem \ref{bau} and ~\cite[Proposition 1.6.8]{a43} to
obtain that for all $\lambda \in {\mathbb C}$ with
$Re\lambda>\beta^2$ and $\widetilde{K_{1}}(\lambda)\neq 0,$ the
following holds
$$
\int \limits^{\infty}_{0}e^{-\lambda t}S(t)xdt=\int
\limits^{\infty}_{0}e^{-\lambda t }\frac{1}{\sqrt{\pi t}}\int
\limits^{\infty}_{0}
e^{-s^{2}/4t}C_{K}(s)xdsdt=\frac{1}{\sqrt{\lambda}}\int
\limits^{\infty}_{0}e^{-\sqrt{\lambda}s}C_{K}(s)xds $$
$$=\frac{1}{\sqrt{\lambda}}\sqrt{\lambda}\tilde{K}(\sqrt{\lambda})(\lambda -A)^{-1}Cx=
\widetilde{K_{1}}(\lambda)(\lambda -A)^{-1}Cx.
$$
As above, one concludes that $(S(t))_{t\geq 0}$ is an
exponentially bounded $K_{1}$-convoluted $C$-semigroup with a
subgenerator $A$. If $Rez>0,$ we define $S(z)x$ in a natural way:
$S(z)x=\frac{1}{\sqrt{\pi z}}\int
\limits^{\infty}_{0}e^{-s^{2}/4z}C_{K}(s)xds.$ Then, $S:\{z\in
{\mathbb C} : Rez>0\} \rightarrow L(E)$ is analytic. Using the
same arguments as in the proof of the Weierstrass formula, see for
instance \cite{a43},  one obtains that for all $\beta \in
(0,\frac{\pi}{2}),$ there exist $M,\ \omega_{1}>0$ such that $
||S(z)||\leq Me^{\omega_{1}|z|},\ z\in \Sigma_{\beta}. $ It
remains to be shown that, for every fixed $\beta \in
(0,\frac{\pi}{2}),$ $\lim_{z\rightarrow 0,\ z\in \Sigma_{\beta}}
S(z)x=0.$ For this, choose an $\omega_{2}>\frac{\omega_{1}}{\cos
\beta}$. Then the function $z\mapsto e^{-\omega_{2}z}S(z)x,\ z\in
\Sigma_{\beta}$ is analytic and satisfies  $\sup_{z\in
\Sigma_{\beta}} ||e^{-\omega_{2}z}S(z)||<\infty.$ Since
$\lim_{t\rightarrow 0+} e^{-\omega_{2}t}S(t)x=0$,
~\cite[Proposition 2.6.3]{a43} implies $\lim_{z\rightarrow 0,\
z\in \Sigma_{\beta}} e^{-\omega_{2}z}S(z)x=0.$ The proof is now
completed.


\section{    Relations to ultradistribution and hyperfunction
sines }\label{bekica}

In this section, we assume $C=I.$ The next assertion clarifies some
properties of generators of local $K$-convoluted cosine functions in
terms of the asymptotic behavior of $\tilde{K}.$

\begin{thm}\label{melj}

\begin{itemize}
\item[(a)]
Suppose that $\Theta$ fulfills $(P2)$ and that $|\Theta(t)|\leq
Me^{\beta t},\ t \geq 0,$ for some $M>0$ and $\beta>0.$ Let $A$ be
the generator of a $K$-cosine function $(C_{K}(t))_{t\in [0,\tau)},$
for some $ \tau\in (0,\infty).$ Further, suppose  that for every
$\varepsilon>0,$ there exist $T_{\varepsilon}>0$ and
$\varepsilon_{0}\in (0,\tau \varepsilon)$ such that
$$
\frac{1}{|\tilde{\Theta}(\lambda)|}\leq
T_{\varepsilon}e^{\varepsilon_{0}|\lambda|},\; \lambda \in
\Omega_{\varepsilon,C_\varepsilon} \cap \{\lambda \in {\mathbb C} :
Re\lambda>\beta\}.
$$
Then, for every $\varepsilon>0,$ there exist positive real numbers
$\overline{C}_{\varepsilon}$ and $\overline{K}_{\varepsilon}$ such
that
$$
 \Omega_{\varepsilon,\overline{C}_{\varepsilon}}^{2}
\subset \rho(A)\mbox{  and } ||R(\lambda^{2}:A)||\leq
\overline{K}_{\varepsilon}e^{\tau \varepsilon |\lambda|},\mbox{
}\lambda \in \Omega_{\varepsilon,\overline{C}_{\varepsilon}}.
$$
\item[(b)] Let $|\Theta(t)|\leq Me^{\beta t},\ t\geq 0,$ for some
$M>0$ and $\beta>0.$ Let $K$ satisfy $(P2).$ Assume that the
restriction of $K$ on $[0,\tau)$ (denoted by the same symbol) is
$\neq 0$ and that $A$ is the generator of a local $K$- convoluted
cosine function on $[0,\tau)$. If there is an $\alpha>0$ with $
\frac{1}{\tilde{|\Theta}(\lambda)|}=O(e^{M(\alpha \lambda)}),
\;|\lambda |\rightarrow \infty ,$ then, for every $\tau_{1} \in (0,
\tau),$ there exist $\beta>0$ and $ C>0$ such that
$$
\Lambda^{2}_{\alpha, \beta,\tau_{1}}\subset \rho(A) \mbox{ and that
} ||R(\lambda:A)||\leq C\frac{e^{M(\alpha
\sqrt{\lambda})}}{1+\sqrt{|\lambda|}},\mbox{ }\lambda \in
\Lambda^{2}_{\alpha, \beta, \tau_{1}}.
$$
\end{itemize}
\end{thm}
{\em Proof:}
(a) Since ${\mathcal A}$ generates a $\Theta$-semigroup
$(S_{\Theta}(t))_{t\in [0,\tau)}$ in $E^{2}$, then we have proved in
\cite{ms} that there are an $\overline{C}_{\varepsilon}>0$ and an $
\overline{K}_{\varepsilon}>0$ so that
$\Omega_{\varepsilon,\overline{C}_{\varepsilon}}\subset
\rho({\mathcal A})$ and $ ||R(\lambda:{\mathcal A})||\leq
\overline{K}_{\varepsilon}e^{\tau \varepsilon |\lambda|},\mbox{
}\lambda \in \Omega_{\varepsilon,\overline{C}_{\varepsilon}}.$  It
follows $\Omega^{2}_{\varepsilon,\overline{C}_{\varepsilon}}\subset
\rho(A)$ and \\ $ ||R(\lambda^{2}:A)||\leq ||R(\lambda:{\mathcal
A})||\leq \overline{K}_{\varepsilon}e^{\tau \varepsilon
|\lambda|},\mbox{   }\lambda \in
\Omega_{\varepsilon,\overline{C}_{\varepsilon}}.$ This finishes the
proof of (a).

(b) We have that ${\mathcal A}$ generates a local $\Theta$-semigroup
on $[0, \tau)$. The prescribed assumption on $\Theta$ and the
arguments of ~\cite[Theorem 1.3.1]{me152} (see also \cite{ms}) imply
that, for every $\tau_{1}\in (0, \tau)$ there exist $\beta>0$ and $
C>0$ such that $\Lambda_{\alpha, \beta, \tau_{1}}\subset
\rho({\mathcal A})$ and that $ ||R(\lambda:{\mathcal A})||\leq
Ce^{M(\alpha \lambda)},\ \lambda \in \Lambda_{\alpha, \beta,
\tau_{1}}. $ Now the proof follows by the standard arguments.

We refer to \cite{ms} for the notion of an ultradistribution
fundamental solution for a closed linear operator $A.$ The notion of
a Fourier hyperfunction fundamental solution for a closed linear
operator $A$ was introduced by Y. Ito in \cite{ito1} while S.
\={O}uchi was the first who introduced the notion of fundamental
solution in the spaces of compactly supported hyperfunctions (cf.
\cite{o192}).

For the sake of simplicity, we use the next definition of
ultradistribution and (Fourier) hyperfunction sines employed by H.
Komatsu in \cite{k92} in the case of an ultradistribution sine.
Similarly, one can introduce and prove the basic characterizations
of tempered ultradistribution sines (cf. \cite{ms}).

\begin{defn} \label{sin}
A closed operator $A$ generates an ultradistribution sine of
$*-$class if there exists an ultradistribution fundamental solution
for the operator ${\mathcal A}.$ A closed operator $A$ generates a
(Fourier) hyperfunction sine if there exists a (Fourier)
hyperfunction fundamental solution for ${\mathcal A}.$
\end{defn}

\begin{rem} \em We will not
 go into  details concerning a relationship between
ultradistribution (hyperfunction) sines and the solvability of
convolution type equations in vector-valued ultradistribution
(hyperfunction) spaces. This can be a matter of further
investigations. In the case of distribution cosine functions, such
an analysis is obtained in \cite{ko99} by passing to the theory of
distribution semigroups (see ~\cite[Theorem 3.10]{ko99}). It is not
so straightforward to link ultradistribution (hyperfunction) sine
generated by $A,$ denoted by $G,$ with ultradistribution fundamental
solution for ${\mathcal A},$ denoted by ${\mathcal G}.$ In
distribution case, we have ${\mathcal G}= \left(
\begin{array}{cc}
G & G^{-1} \\
G^{\prime }-\delta & G
\end{array}
\right)$ (see \cite{ko99} for more details). The main problem in
transferring ~\cite[Theorem 3.10(i)]{ko99} to ultradistribution and
hyperfunction sines is the appearing of $G^{-1}$ in the
representation formula for ${\mathcal G}.$ Furthermore, relations
between (almost-)distribution cosine functions and cosine
convolution products have been recently analyzed in \cite{mp1} and
\cite{miki}. It is not clear how to obtain the corresponding results
in the case of ultradistribution and hyperfunction sines.
\end{rem}

Spectral properties of operators generating ultradistribution and
(Fourier) hyperfunction sines are given in the following remark.

\begin{rem} \label{uporedno}
\em
 1 (\cite{k92}, \cite{ku113}). A closed linear operator $A$
generates an ultradistribution sine of $*-$class iff there exists a
domain of the form $\Omega^*$ such that:
\begin{equation} \label{relac1}
\Omega^{*,2}\subset \rho(A) \mbox{   and }
\end{equation}
\begin{equation} \label{relac2}
||R(\lambda^{2}:A)||\leq Ce^{M(k|\lambda|)},\mbox{   } \lambda\in
\Omega^{M_p}_{k,C},
\end{equation}
for some $k>0$ and $ C>0$ in $(M_p)$-case, resp.,
\begin{equation} \label{relac3}
 ||R(\lambda^{2}:A)||\leq C_ke^{M(k|\lambda|)},\mbox{   } \lambda\in \Omega^{M_p}_{k,C_k},
\end{equation}
for every $k>0$ and the corresponding $C_k>0$ in $\{M_p\}$-case.
\vspace{0.1cm}

2 (\cite{ito1}). A closed linear operator $A$ generates a Fourier
hyperfunction sine iff
\begin{equation}\label{itica}
{\mathbb C} \mbox{ }\setminus \mbox{ } (-\infty,0] \subset
\rho(A)\mbox{ and if, for every }\varepsilon>0 \mbox{ and }\sigma>0,
\end{equation}
\begin{equation}\label{temperica}
\mbox{ there is a }C_{\varepsilon,\sigma}>0\mbox{ with }
||R(\lambda^{2}: A)||\leq C_{\varepsilon,\sigma}e^{\varepsilon
|\lambda|},\ Re\lambda>\sigma.
\end{equation}

3 (\cite{o192}). A closed linear operator $A$ generates a
hyperfunction sine iff for every $\varepsilon>0,$ there exist
constants $C_{\varepsilon}>0$ and $ K_{\varepsilon}>0$ satisfying
\begin{equation}\label{boli}
\Omega^{2}_{\varepsilon,C_\varepsilon}\subset\rho(A) \mbox{ and }
||R(\lambda^{2}:A)||\leq K_{\varepsilon}e^{\varepsilon |\lambda|
},\mbox{  }\lambda \in \Omega_{\varepsilon,C_\varepsilon}.
\end{equation}

We refer to \cite{k92} for the spectral properties of operators
generating Laplace hyperfunction semigroups and sines.
\end{rem}

\begin{thm}\label{relacije}

\begin{itemize}
\item[(a)] Let $A$ generate an ultradistribution sine of
$(M_{p})-$class, resp., $\{ M_{p} \}-$class. Then, for every $\theta
\in [0,\frac{\pi}{2}),$ there exists an ultradistribution
fundamental solution of $(M_{p}^{2})-$class, resp., $\{ M_{p}^{2}
\}-$class for $e^{\pm i\theta}A.$
\item[(b)]  Let $A$ generate a
hyperfunction sine. Then, for every $\theta \in [0,\frac{\pi}{2}),$
there exists an ultradistribution fundamental solution of
$\{p!^2\}-$class for $e^{\pm i\theta}A.$
\end{itemize}
\end{thm}
{\em Proof:}
We will prove only (a) since the same arguments work for (b). Fix a $\theta \in
[0,\frac{\pi}{2}).$ If $\overline{M}$ denotes the associated
function of $(M_{p}^{2}),$ then, for every $k>0,$
$M(k\sqrt{t})=\sup\{\ln\frac{M_0k^{p}\sqrt{t^p}}{M_p} : p\in
{{\mathbb N}_{0}} \} =\frac{1}{2}\sup\{\ln\frac{M_0^2
k^{2p}t^p}{M^2_p} : p\in {\mathbb N}_0\}=\frac{1}{2}
\overline{M}(k^{2}t),$ $t\geq 0.$ We have already noted that $A$
generates an ultradistribution sine of $(M_p)-$class, resp.,
$\{M_{p}\}-$class if and only if there exists a domain of the form
$\Omega^*$ such that (\ref{relac1}) and (\ref{relac2}), resp.,
(\ref{relac1}) and (\ref{relac3}) are fulfilled. Since
$$
\lim \limits_{|\lambda|\rightarrow \infty,\ \lambda \in \partial
(\Omega^{ M_{p}}_{k,C})}\cos(\arg(\lambda))=\lim
\limits_{|\lambda|\rightarrow
\infty}\frac{M(k|\lambda|)+C}{|\lambda|}=0,
$$
we obtain $|\arg \lambda|\rightarrow \frac{\pi}{2},\; | \lambda |
\rightarrow \infty,\ \lambda \in
\partial (\Omega^{M_{p}}_{k,C}),$ and therefore, $|\arg
\lambda|\rightarrow \pi,\; | \lambda | \rightarrow \infty,\ \lambda
\in
\partial ((\Omega^{M_{p}}_{k,C})^2).$ The same estimate holds in the Roumieu case. Hence, there
exists a suitable $\omega>0$ with $\rho(A)\supset\Omega^{*,2}
\supset \omega +\Sigma_{\pi/2+\theta}.$ Further,\\ $
||R(\lambda:A)||\leq Ce^{M(k\sqrt{|\lambda|})},\mbox{   } \lambda\in
(\Omega^{M_p}_{k,C})^{2},$ for some $k>0$ and $ C>0$ in
$(M_p^2)-$case, resp., $ ||R(\lambda:A)||\leq
C_ke^{M(k\sqrt{|\lambda|})},\mbox{   } \lambda\in
(\Omega^{M_p}_{k,C_k})^{2},$ for every $k>0$ and the corresponding
$C_k>0,$ in $\{M_p^2\}-$case. The analysis given in the first part
of the proof shows that, for every $\theta \in [0,\frac{\pi}{2}),$
we have $\{ z\in {\mathbb C} : Rez>\omega \} \subset \rho(e^{\pm
i\theta}A),$ and that $||R(\lambda : e^{\pm i\theta}A)||=||R(\lambda
e^{\mp i\theta} : A)||\leq
Ce^{M(k\sqrt{|\lambda|})}=Ce^{\frac{1}{2}\overline{M}(k^{2}|\lambda|)},\mbox{
} Re\lambda>\omega,$ in $(M_p^2)-$case. The similar estimate holds
in the Roumieu case. Now the proof  finishes an application of
arguments given in ~\cite{ms}.

Next, we relate ultradistribution and hyperfunction sines to
analytic convoluted semigroups. Recall, the function
$K_{1/2}(t)=\frac{1}{2\sqrt{\pi t^{3}}}e^{-\frac{1}{4t}}$, $t> 0$ is
bounded and smooth. Furthermore, $\tilde{K}(\lambda
)=e^{-\sqrt{\lambda }},$ $Re\lambda>0,$ where $\sqrt{1}=1.$

\begin{thm}\label{anaveze}
Suppose that $A$ generates a hyperfunction sine. Then $A$ generates
an exponentially bounded, analytic $K_{1/2}$-semigroup of angle
$\frac{\pi}{2}$.
\end{thm}
{\em Proof:}
 Clearly, it is enough to show that, for every $\varepsilon
\in (0, \frac{1}{\sqrt{2}})$, $A$ generates an exponentially bounded
analytic $K_{1/2}$-semigroup of angle
$\alpha:=2\arccos{\varepsilon}-\frac{\pi}{2}$. So let $\varepsilon
\in (0,\frac{1}{\sqrt{2}})$ be fixed. Then there exist
$C_{\varepsilon}>0$ and $ K_{\varepsilon}>0$ such that\\
$\Omega^{2}_{\varepsilon, C_{\varepsilon}}\subset \rho(A)$ and that
$ ||R(\lambda:A)||\leq K_{\varepsilon}e^{\varepsilon
\sqrt{|\lambda|}},\ \lambda \in
\Omega_{\varepsilon,C_\varepsilon}^{2}.$ Since
$$
\partial
(\Omega_{\varepsilon,C_\varepsilon})=\{ re^{i\theta} : r>0,\ \theta
\in (0,\frac{\pi}{2}),\ r\cos \theta=\varepsilon r+C_{\varepsilon}
\},
$$
one can conclude that $ |\arg{\lambda}|\rightarrow
2\arccos{\varepsilon},\ |\lambda|\rightarrow \infty,\ \lambda \in
\partial({\Omega_{\varepsilon,C_\varepsilon}^{2}}).$
This implies that, for a sufficiently large $\omega \in (0,\infty),$
$ \omega+\Sigma_{\frac{\pi}{2}+\alpha}\subset \rho(A).$ Furthermore,
$\lim_{\lambda \rightarrow +\infty}\lambda
\widetilde{K}_{1/2}(\lambda)R(\lambda:A)=0$ and $\tilde{K}_{1/2}$
can be analytically extended to the function $g:
\omega+\Sigma_{\frac{\pi}{2}+\alpha} \to {\mathbb C}$,
$g(\lambda)=e^{-\sqrt{\lambda}},\ \lambda \in
\omega+\Sigma_{\frac{\pi}{2}+\alpha}$. Fix a $\gamma_{1}\in
(0,\alpha).$ Then it is straightforward to see that
$\cos(\frac{\pi}{4}+\frac{\gamma_{1}}{2})>\varepsilon$ and that
$$
||(\lambda-\omega)g(\lambda)R(\lambda:A)||\leq
K_{\varepsilon}(|\lambda|+\omega)e^{\varepsilon
\sqrt{|\lambda|}}e^{-\sqrt{|\lambda|}\cos(\arg(\lambda) /2)} $$
$$
\leq
K_{\varepsilon}(|\lambda|+\omega)e^{\sqrt{|\lambda|}(\varepsilon-\cos(\frac{\pi}{4}+\frac{\gamma_{1}}{2}))},\
\lambda \in \omega+\Sigma_{\gamma_{1}+\frac{\pi}{2}}. $$
Theorem
\ref{anali} ends the proof.

The similar assertion holds for ultradistribution sines. For the
sake of brevity, in the next theorem, we consider only the case when
$(M_{p})$ is a Gevrey type sequence: $(p!^{s}),\ (p^{ps})$ or
$(\Gamma(1+ps)),\ s>1.$ Then we know that, for every $s>1,$ there
exists an appropriate $C_{s}^{\prime}>0$ so that $M(t)\sim
C_{s}^{\prime}t^{\frac{1}{s}},\ t\rightarrow +\infty.$

\begin{thm}\label{anault}
Suppose that $A$ generates an ultradistribution sine of the
Beurling, resp.,  Roumieu class. Then $A$ generates an exponentially
bounded, analytic $K_{\delta}$-semigroup of angle $\frac{\pi}{2}$,
for all $\delta \in (\frac{1}{2s}, \frac{1}{2}),$ resp., for all
$\delta \in [\frac{1}{2s}, \frac{1}{2}).$
\end{thm}
{\em Proof:}
 We prove the assertion in the Roumieu case since the proof in
the Beurling case can be derived similarly. Let us fix some $\gamma
\in (0, \frac{\pi}{2})$ and $\delta \in [\frac{1}{2s},
\frac{1}{2}).$ We know that, for every $k>0$ and a corresponding
$C_k>0:$
$$
\{\lambda^{2}:\lambda \in \Omega^{M_p}_{k,C_k} \}\subset \rho(A)
\mbox{   and } \; ||R(\lambda^{2}:A)||\leq
C_ke^{M(k|\lambda|)},\mbox{   } \lambda\in \Omega^{M_p}_{k,C_k}.
$$
Since $M(|\lambda|)\leq C_{s}|\lambda|^{1/s},\ |\lambda| \geq 0$ for
some $C_{s}>0,$  one obtains
$$\{\lambda^{2} : Re\lambda \geq C_{s}(k|\lambda|)^{1/s}+{C}_k\}
\subset \rho(A), \mbox{ i.e.,}$$
$$\{ r^{2}e^{2i\theta}: r>0,\
|\theta|<\frac{\pi}{2},\ r\cos \theta \geq C_{s}k^{1/s}r^{1/s}+{C}_k
\} \subset \rho(A).$$ Denote $\Gamma =\{re^{i\theta} : r\cos \theta
= C_{s}k^{1/s}r^{1/s}+{C}_k\}.$ Then $\lim_{|\lambda|\rightarrow
\infty, \lambda \in \Gamma} |\arg (\lambda)|=\frac{\pi}{2}.$
Therefore, there are an $\omega_{\gamma}>0$ and a suitable
$\overline{C}_{k}>0$ so that $\omega_{\gamma} +
\Sigma_{\frac{\pi}{2} +\gamma}\subset \rho(A)$ and that
$$
||R(\lambda : A)||\leq \overline{C}_{k} e^{M(k\sqrt{|\lambda|})}\leq
\overline{C}_{k}
 e^{C_{s}k^{1/s}|\lambda|^{1/2s}},\ \lambda \in
\omega_{\gamma} + \Sigma_{\frac{\pi}{2} +\gamma}.
$$
Clearly, the function $g: \omega_{\gamma} + \Sigma_{\frac{\pi}{2}
+\gamma} \rightarrow {\mathbb C}$, $g(\lambda)=e^{-\lambda
^{\delta}},$ $1^{\delta}=1$ is analytic. Furthermore,
$$ ||g(\lambda)R(\lambda
: A)||\leq \frac{\overline{C}_{k}}{|e^{\lambda
^{\delta}}|}e^{C_{s}k^{1/s}|\lambda|^{1/2s}}$$
$$=
\overline{C}_{k}\exp{(-|\lambda|^{\delta}\cos(\delta \arg (\lambda))
+C_{s}k^{1/s}|\lambda|^{1/2s})}$$
$$ \leq \overline{C}_{k} e^{C_{s}k^{1/s}|\lambda|^{1/2s}-\cos(\pi
\delta) |\lambda|^{\delta}},\ \lambda \in \omega_{\gamma} +
\Sigma_{\frac{\pi}{2} +\gamma}. $$
 Our choice of $\delta$, the
fact that a number $k>0$ can be chosen arbitrarily in the Roumieu
case and Theorem \ref{anali}, imply that $A$ generates an
exponentially bounded, analytic $K_{\delta}$-semigroup of angle
$\gamma$. This ends the proof.

Motivated by \cite{keya} and \cite{mp1}, up to the end of this
section, we discuss relations between (local) integrated cosine
functions and ultradistribution semigroups (sines). Recall, a closed
linear operator $A$ generates a local integrated cosine function if
and only if there exist $\alpha,\ \beta,\ M>0$ and $n\in {\mathbb
N}$ so that $E^{2}(\alpha, \beta) \subset \rho(A)$ and $||R(\lambda
: A)||\leq M(1+|\lambda|)^{n},\ \lambda \in E^{2}(\alpha,\beta)$
(see \cite{ko99}).

\begin{rem}\label{kos}
\em
 We recall that V. Keyantuo proved in ~\cite[Theorem 3.1]{keya}
that if a densely defined operator $A$ generates an exponentially
bounded $\alpha$-times integrated cosine function for some $\alpha
\geq 0$ (this means that $A$ is densely defined and that $A$
generates an exponential distribution cosine function of
\cite{ko99}), then $\pm iA$ generate ultradistribution semigroups of
\cite{cha}. Furthermore, the proof of ~\cite[Theorem 3.1]{keya} and
the assertion of ~\cite[Proposition 2.6]{ci1} imply that $\pm iA$
generate regular $(M_{p})-$ultradistribution semigroups in the sense
of ~\cite[Definition 2.1]{ci1}, where $M_{p}=p!^{s},$ $s\in (1,2);$
see also Example \ref{Beals} given below. In general, the assertion
of ~\cite[Theorem 3.1]{keya} does not hold if $s>2$ and, in the case
$s=2$, this assertion remains true only in the Beurling case, see
~\cite[Section 4]{mp1}.

In \cite{mp1}, the next extension of ~\cite[Theorem 3.1]{keya} has
been recently showed: If $A$ is the generator of a (local)
$\alpha$-times integrated cosine function, for an $\alpha>0,$ then
$\pm iA$ are generators of ultradistribution semigroups of
$\ast-$class, where $(M_{p})$ is a Gevrey type sequence with $s\in
(1,2).$ Then it can be easily seen that $-A^{2}$ generates an
ultradistribution sine of $({M_{p}})-$class, resp., $\{ M_{p}
\}-$class.

We want to notice that, in general, $-A^{2}$ does not generate a
local integrated cosine function even if $A$ is the densely defined
generator of an exponentially bounded, integrated cosine function.
In order to illustrate this fact, we choose $E:=L^{p}({{\mathbb
R}}),\ 1\leq p<\infty$ and put $m(x):=(1-\frac{x^{2}}{4})+ix,\ x\in
{\mathbb R}.$ Define a closed linear operator $A$ on $E$ by:
$$
Af(x)=m(x)f(x),\ x\in {\mathbb R},\ D(A):=\{ f\in E : mf\in E \}.
$$
Then it is proved in ~\cite[Example 4.4]{mp1} that $A$ generates a
dense, exponential distribution cosine function. We have
$$
\sigma (-A^{2})=\{
(\frac{x^{4}}{16}+\frac{3}{2}x^{2}-1)+2ix(\frac{x^{2}}{4}-1) : x\in
{\mathbb R} \}
$$
and a simple analysis shows that there do not exist $\alpha>0$ and $
\beta>0$ with \\ $E^{2}(\alpha,\beta) \cap \sigma(-A^{2}) \neq
\emptyset$ (see ~\cite[Section 4]{mp1} for more details). Hence,
$-A^{2}$ does not generate a local $\alpha$-times integrated cosine
function, for any $\alpha>0.$
\end{rem}

The next theorem improves ~\cite[Theorem 3.1]{keya} in a different
direction.

\begin{thm}\label{dedicated}
Let $K$ satisfy (P1) and $(P2)$ and let $A$ be the generator of an
exponentially bounded, $K$-cosine function $(C_{K}(t))_{t\geq 0 }.$
Suppose $M_{p}=p!^{s},$ for some $s\in (1,2).$ If there exist $k>0$
and $\overline{C}>0,$ in the Beurling case, resp., if for every
$k>0$ there exists an appropriate $\overline{C_{k}}>0,$ in the
Roumieu case, such that
\begin{equation}\label{mance1}
\frac{1}{|\tilde{K}(\lambda)|}=O(e^{M(k|\lambda|)}),\ Re\lambda \geq
\overline C,\mbox{  resp.,  }Re\lambda \geq \overline{C_{k}},
\end{equation}
then there exist ultradistribution fundamental solutions of
$\ast-$class for $\pm iA.$
\end{thm}
{\em Proof:}
We prove the assertion in the Roumieu case. To do this, fix $k>0.$
We know that there exist $a> \overline{C_{k}},$ $l>0$ and $L>0$ with
$e^{l|\lambda|^{1/s}}\leq e^{M(|\lambda|)}\leq
e^{L|\lambda|^{1/s}},$ $\; \lambda \in {\mathbb C},\ |\lambda|\geq
a.$ Theorem \ref{bau} and the assumption (\ref{mance1}) imply that
there exists an $\omega >$max(abs$(K),a)$ so that $ ||R(\lambda :
A)||=O(e^{M(k\sqrt{|\lambda|})}),\ \lambda \in \Pi_{\omega}.$
Since $\partial \Pi_{\omega}=\{x+iy \in {\mathbb C} : x = \omega ^{2}-\frac{%
y^{2}}{4\omega ^{2}}\},$ we have
\begin{equation}\label{bokacc}
\{z\in \mathbb{C} : Imz\geq (\omega
+1)^{2}-\frac{(Rez)^{2}}{4(\omega +1)^{2}}\} \subset \rho(iA)
\end{equation}
and for such $\lambda$'s: $||R(\lambda : iA
)||=O(e^{M(k\sqrt{|\lambda|})}).$ The choice of $s$ implies that
there exists a sufficiently large $\beta>0$ with
$$\frac{x^{2}}{4(\omega +1)^{2}}-\frac{x^{s}}{kl^{s}}\geq (\omega
+1)^{2},\ x\geq \beta.$$ Put now $C_{k}=\max (\frac{a}{k}, \beta).$
Suppose
$$z=x+iy\in \Omega_{k,{C_{k}}}^{M_{p}},\mbox{  i.e., } Rez\geq
M(k|z|)+C_{k}.$$ Then $y^{2}\leq
(\frac{x-C_{k}}{lk^{1/s}})^{2s}-x^{2}.$ According to the choice of
$C_{k},$ one obtains
$$
y+\frac{x^{2}}{4(\omega +1)^{2}}\geq \frac{x^{2}}{4(\omega
+1)^{2}}-\sqrt{(\frac{x-C_{k}}{lk^{1/s}})^{2s}-x^{2}}$$
$$
\geq \frac{x^{2}}{4(\omega +1)^{2}}-(\frac{x-C_{k}}{lk^{1/s}})^{s}
\geq \frac{x^{2}}{4(\omega +1)^{2}}-(\frac{x}{lk^{1/s}})^{s}\geq
(\omega +1)^{2}.
$$
Due to (\ref{bokacc}), $z\in \rho(iA).$ We know $||R(\lambda : i
A)||=O( e^{M(k\sqrt{|\lambda|})}),\ \lambda \in
\Omega^{M_p}_{k,C_{k}}$ and this proves the claimed assertion for
$iA.$ The same arguments work for $-iA.$

Suppose $M_{p}=p!^{2}$ in the formulation of Theorem
\ref{dedicated}. As before, this theorem remains
true only in the Beurling case. 

\section{  Examples and applications }

\begin{example}\label{biha}
\em
 Let $A := -\Delta $ on $E := L^{2}[0,\pi ]$ with the Dirichlet
boundary conditions (see, for example, ~\cite[Section 7.2]{a43}).
Motivated by the paper of B. B\" aumer \cite{b40} we have proved in
\cite{ms} that there exists an exponentially bounded kernel $K\in
C([0,\infty))$ so that $A$ generates a $K$-semigroup
$(S_{K}(t))_{t\geq 0}$ with $||S_{K}(t)||=O(1+t^{3})$. Suppose that
$|K(t)|\leq Me^{\beta t},\ t\geq 0,$ for some $M>0$ and $\beta>0.$
Moreover, $-A$ also generates an exponentially bounded $K$-semigroup
$(V_{K}(t))_{t\geq 0}$ in $E$ since it is the generator of an
analytic $C_{0}$-semigroup of angle $\frac{\pi}{2}.$ Then
Proposition \ref{veza} implies that the biharmonic operator
$\Delta^{2},$ endowed with the corresponding boundary conditions,
generates an exponentially bounded, $K$-cosine function
$(C_{K}(t))_{t\geq 0},$ where
$C_{K}(t):=\frac{1}{2}(S_{K}(t)+V_{K}(t)),\ t\geq 0.$ Put
$K_{1}(t):=\int \limits^{t}_{0}K(s)ds,\ t\geq 0.$ Clearly,
$|K_{1}(t)|\leq Mte^{\beta t},\ t\geq 0,$ and $\Delta^{2}$ generates
an exponentially bounded, $K_{1}$-cosine function
$(C_{K_{1}}(t))_{t\geq 0},$ where $C_{K_{1}}(t)x=\int
\limits^{t}_{0}C_{K}(s)xds,\ x\in E,\ t\geq 0.$ This implies that
$\Delta^{2}$ generates an exponentially bounded, analytic
$K_{2}$-semigroup of angle $\frac{\pi}{2}$, where $K_{2}(t):=\int
\limits^{\infty}_{0}\frac{se^{-s^{2}/4t}}{2\sqrt{\pi}
t^{3/2}}{K_{1}}(s)ds,\ t>0.$ Note that we have integrated once the
function $K$ in order to prove that $K_{2}$ is exponentially
bounded. This is valid, since for every $t>0:$
$$
|K_{2}(t)|\leq M\int
\limits^{\infty}_{0}\frac{se^{-s^{2}/4t}}{2\sqrt{\pi} t^{3/2}}
se^{\beta s}ds = \frac{M}{2\sqrt{\pi}}\int
\limits^{\infty}_{0}r^{2}e^{\beta r \sqrt{t}-\frac{r^{2}}{4}}dr
$$
$$
= \frac{M}{2\sqrt{\pi}}\int
\limits^{\infty}_{0}r^{2}e^{\beta^{2}t-(\frac{r}{2}-\beta
\sqrt{t})^{2}}dr
$$
$$
=\frac{M}{2\sqrt{\pi}}e^{\beta^{2}t}\int
\limits^{\infty}_{0}r^{2}e^{-(\frac{r}{2}-\beta
\sqrt{t})^{2}}dr=\frac{M}{2\sqrt{\pi}}e^{\beta^{2}t}\int
\limits^{\infty}_{-\beta \sqrt{t}}8(v^{2}+2v\beta
\sqrt{t}+\beta^{2}t)e^{-v^{2}}dv
$$
$$
\leq \frac{4M}{\sqrt{\pi}}e^{\beta^{2}t}[\int
\limits^{\infty}_{-\infty}v^{2}e^{-v^{2}}dv+2\beta \sqrt{t}\int
\limits^{\infty}_{0}ve^{-v^{2}}dv+\beta^{2}t\int
\limits^{\infty}_{-\infty}e^{-v^{2}}dv]\leq
\overline{M}e^{(\beta^{2}+1)t},
$$
for a suitable $\overline{M}>0.$ Furthermore, $K_{2}$ is a kernel
since
$$
\limsup_{\lambda \rightarrow \infty}
\frac{\ln{|\widetilde{K_{2}}(\lambda)|}}{\lambda}=\limsup_{\lambda
\rightarrow \infty}
\frac{\ln{|\widetilde{K_{1}}(\sqrt{\lambda})|}}{\lambda}=0.
$$
On the other side, $\Delta^{2}$ cannot be the generator of a (local)
integrated $\alpha$-times semigroup, $\alpha \geq 0$, since the
resolvent set of $\Delta^{2}$ does not contain any ray $(\omega,
\infty)$. For the same reasons, $\Delta^{2}$ does not generate a
hyperfunction (ultradistribution) sine. Hence, in the analysis of
$\Delta^{2}$ and $-\Delta,$ we do not need any $C,$ but the use of
regularized operator families enables several advantages which
hardly can be considered by the use of asymptotic Laplace transform
techniques. More generally, suppose $n\in {\mathbb N}$. Since
$\Delta=-A$ generates a cosine function (see  ~\cite[Example 7.2.1,
p. 418]{a43}), one can employ a result of J. A. Goldstein proved in
\cite{goja} (see also ~\cite[p. 215]{l1}), in order to see that
$-\Delta^{2n}$ generates an analytic $C_{0}$-semigroup of angle
$\frac{\pi}{2}.$ Hence, an application of ~\cite[Theorem 8.2]{l1}
shows that there exists an injective operator $C_{n}\in
L(L^{2}[0,\pi])$ so that $\Delta^{2n}$ generates an entire
$C_{n}$-group. Further on, one can apply Proposition \ref{veza} in
order to see that the polyharmonic operator $\Delta^{4}$ generates
an exponentially bounded, $K_{2}$-convoluted cosine function. Put
$\overline{K_{3}}(t):=\int \limits^{t}_{0}K_{2}(s)ds,\ t\geq 0.$
Then $K_{3}$ is a kernel and we have $|\overline{K_{3}}(t)|\leq
\overline{M}te^{(\beta^{2}+1)t},\ t\geq 0.$ Clearly, $\Delta^{4}$
generates an exponentially bounded, $\overline{K_{3}}$-cosine
function. Then Theorem \ref{mance} can be applied again in order to
see that $\Delta^{4}$ generates an exponentially bounded, analytic
$K_{3}$-semigroup of angle $\frac{\pi}{2},$ where $K_{3}(t):=\int
\limits^{\infty}_{0}\frac{se^{-s^{2}/4t}}{2\sqrt{\pi}
t^{3/2}}\overline{K_{3}}(s)ds,\ t>0.$ Similarly as above, we have
that $K_{3}$ is an exponentially bounded kernel. Continuing this
procedure leads us to the fact, mentioned already in the abstract
and the introduction of the paper, that there exist exponentially
bounded kernels $K_{n}$ and $K_{n+1}$ such that $\Delta^{2^{n}}$
generates an exponentially bounded, $K_{n}$-convoluted cosine
function, and in the meantime, an exponentially bounded, analytic
$K_{n+1}$-convoluted semigroup of angle $\frac{\pi}{2}.$ Note that
this procedure can be done only with the loos of regularity, since
we must apply Theorem \ref{mance} (see also~\cite[Proposition
1.6.8]{a43}). At the end of this discussion, note that it is not
clear whether there exists a kernel $\overline{K_{n}}$ such that
$\Delta^{2n}$ generates an exponentially bounded,
$\overline{K_{n}}$-convoluted cosine function.

Suppose now that $A$ is a self-adjoint operator in a Hilbert space
$H$ and that $A$ has a discrete spectrum $(\lambda_{n})_{n\in
{\mathbb N}},$ where we write the eigenvalues in increasing order
and repeat them according to multiplicity. Suppose that $Re
\lambda_{n}>0,\ n\geq n_{0}$ and $m$ is a natural number which is
greater than any multiplicity of $\lambda_{n},\ n\geq n_{0}.$ If
$$
\sum_{n\geq
n_{0}}(1-\frac{|\sqrt{\lambda_{n}}-1|}{\sqrt{\lambda_{n}}+1}
)<\infty,
$$
then, according to ~\cite[Theorem 1.11.1]{a43}, there exists an
exponentially bounded function $K$ such that
$\tilde{K}(\sqrt{\lambda_{n}})=0$, $n\geq n_{0}$. This implies that
the function $\lambda \mapsto \widetilde{K^{\ast
m}}(\lambda)R(\lambda^{2}:A)$ can be analytically extended on a
right half plane, where $K^{\ast m}$ denotes the $m$-th convolution
power of $K$. If, additionally,
$$
||\widetilde{K^{\ast m}}(\lambda)R(\lambda^{2}:A)||\leq
M|\lambda|^{-3},\ Re\lambda>\omega \ (\geq 0),\ \lambda \neq
\sqrt{\lambda_{n}},\ n\geq n_{0},
$$
for a suitable $M>0$, then $A$ generates an exponentially bounded
$K^{\ast m}$-cosine function. The main problem is to construct a
kernel $K$ which fulfills the previous estimate. It is also
evident that this procedure cannot be done if
$(\sqrt{\lambda_{n}})_{n\geq n_{0}}$ is a uniqueness sequence, see
for instance \cite{a43} and \cite{b40}. Therefore, the theory of
convoluted cosine functions cannot be applied if $\lambda_{n} \sim
n^{2s},\ n\rightarrow +\infty,$ for some $s\in (0,1],$ and it, in
turn, implies that the operator $-\Delta,$ considered in the first
part of this example, cannot be the generator of any exponentially
bounded, convoluted cosine function. Finally, we refer to
~\cite[Chapter 2]{me152} for the notion and basic properties of
the spaces of "new distributions", for the treatment of such kind
of problems within the theory of generalized functions.
\end{example}

\begin{example}
\em
 Let ${\mathbb C}_{+}=\{z\in {\mathbb C} : Imz>0\}$ and $1\leq
p<\infty.$ Suppose
that $E:=H^{p}({\mathbb C}_{+}).$ 
Recall that  R. Beals constructed in the proof of ~\cite[Theorem
2']{b41} an analytic function $a_{1} : {\mathbb C}_{+} \to \{ z\in {
{\mathbb C} } : |z| \geq 1 \}$ with the property that, for every
$\varepsilon
>0,$ there exists a region of the form $\Omega_{\varepsilon,C_\varepsilon}$ satisfying
$a_{1}({\mathbb C}_{+}) \cap \Omega_{\varepsilon,C_\varepsilon}
=\emptyset .$ Let $B=a_{1}^{2}$. Then $B$ is a holomorphic function
on ${\mathbb C}_{+}$ and for all $\varepsilon>0$ there exist
$C_{\varepsilon}>0$ and $ K_{\varepsilon}>0$ so that $B({\mathbb
C}_{+})\subset (\Omega^{2}_{\varepsilon,C_\varepsilon})^{c}.$ Define
$$
(AF)(z):=B(z)F(z),\;Imz>0,\mbox{ }D(A):=\{F\in H^{p}({\mathbb C}_{+}):%
AF\in H^{p}({\mathbb C}_{+})\}.
$$
Let $\varepsilon \in (0,1)$ be fixed. Choose an $\varepsilon_{1}\in
(0, \varepsilon)$ such that $B({\mathbb C}_{+})\subset
(\Omega^{2}_{\varepsilon_{1},C_{\varepsilon_1}})^{c}$. Clearly,
$\lim_{|\lambda| \rightarrow \infty, \lambda \in \partial
(\Omega_{\varepsilon_{1},C_{\varepsilon_1}})} |\arg \lambda|=\arccos
\varepsilon_{1}$ and there exists a sufficiently large
$\overline{C}_{\varepsilon}>0$ such that
$\Omega_{\varepsilon,\overline{C}_\varepsilon}=\{\lambda :\; Re
\lambda \geq \varepsilon |\lambda|+ \overline{C}_{\varepsilon}\}
\subset \Omega_{\varepsilon_{1},C_{\varepsilon_1}}$ and that the
distance $d:=$dist$ (\partial (\Omega_{{ \varepsilon
}_{1},C_{\varepsilon_1}}),
\partial (\Omega_{\varepsilon,\overline{C}_\varepsilon}))>0$. This
implies: \\
$\Omega^{2}_{\varepsilon,\overline{C}_\varepsilon}\subset \rho(A)$
and $ ||R(\lambda:A)||\leq d^{-2},\mbox{   }\lambda \in \Omega^{2}_{
\varepsilon,\overline{C}_\varepsilon}. $ Therefore, $A$ generates a
hyperfunction sine, and it can be easily seen that $A$ does not
generate an ultradistribution sine of $\ast-$class.
\end{example}

\begin{example}\label{Beals}
\em
 Let $E=L^{p}({\mathbb R}),\ 1\leq p \leq \infty$. Consider the
next multiplication operator with the maximal domain in $E:$
$$
Af(x):=(x+ix^{2})^{2}f(x),\ x \in {\mathbb R},\ f\in E.
$$
It is clear that $A$ is dense and stationary dense if $1\leq p
<\infty$, but $A$ is not the generator of any (local) integrated
cosine function, $1\leq p \leq \infty$. Moreover, if $p=\infty$,
then $A$ is not stationary dense since, for example, the function
$x\mapsto \frac{1}{x^{2n}+1}$ belongs to $D(A^{n}) \setminus
\overline{D(A^{n+1})},\; n\in {\mathbb N}.$ Further, one can easily
verify that $A$ generates an ultradistribution sine of $\ast-$class,
if $M_{p}=p!^{s},\ s \in (1,2).$ If $M_{p}=p!^{2},$ then an analysis
given in ~\cite[Example 4.4]{mp1} shows that $A$ does not generate
an ultradistribution sine of the Roumieu class and that $A$
generates an ultradistribution sine of the Beurling class. Suppose
now $M_{p}=p!^{s},$ for some $ s\in (1,2),$ and put
$\delta=\frac{1}{s}.$ Then $A$ generates a global (non-exponentially
bounded) $K_{\delta}$-cosine function since, for every $\tau \in
(0,\infty),$ $A$ generates a $K_{\delta}$-cosine function on
$[0,\tau).$ Indeed, suppose $M(\lambda)\leq C_{s}|\lambda|^{1/s},\
\lambda \in {\mathbb C}.$ Fix  $\tau \in (0,\infty)$ and choose an
$\alpha
>0$ with $\tau \leq \frac{\cos(\frac{\delta \pi
}{2})}{C_{s}\alpha^{1/s}}.$ The choice of $\alpha$ implies that
there exists a sufficiently large $\beta>0$ such that
$\Lambda_{\alpha,\beta,1}^{2}\subset \rho(A)$ and that the resolvent
of $A$ is bounded on $\Lambda_{\alpha,\beta,1}^{2}.$ Put
$\Gamma:=\partial (\Lambda_{\alpha,\beta,1}).$ We assume that
$\Gamma$ is upwards oriented. Define
$$
C_{\delta}(t)f(x):=\frac{1}{2\pi i}\int
\limits_{\Gamma}\frac{\lambda e^{\lambda
t-\lambda^{\delta}}}{\lambda^{2}-(x+ix^{2})^{2}}d\lambda f(x),\ f\in
E,\ x\in {\mathbb R},\ t\in [0,\frac{\cos(\frac{\delta \pi
}{2})}{C_{s}\alpha^{1/s}}).
$$
Note that the above integral is convergent since
$$|e^{-\lambda^{\delta}}|\leq e^{-\cos(\frac{\delta
\pi }{2})|\lambda|^{\delta}},\ Re\lambda>0\mbox{ and}
$$
$$|e^{\lambda
t-\lambda^{\delta}}|\leq e^{\beta t}e^{M(\alpha
\lambda)t-\cos(\frac{\delta \pi}{2})|\lambda|^{\delta}}\leq e^{\beta
t} e^{C_{s}\alpha^{1/s}|\lambda|^{1/s}t-\cos(\frac{\delta
\pi}{2})|\lambda|^{\delta}},\ \lambda \in \Gamma.
$$
It is straightforward to check that $(C_{\delta}(t))_{t\in
[0,\tau)}$ is a $K_{\delta}$-cosine function generated by $A.$ At
the end, we point out that there exists an appropriate $\tau_{0} \in
(0,\infty)$ such that $A$ generates a local $K_{1/2}$-cosine
function on $[0,\tau_{0}).$
\end{example}

Many other examples of differential operators, acting on
$L^{2}({\mathbb R}^{n}),\ n\in {\mathbb N},$  which generate
ultradistribution and hyperfunction sines (semigroups) can be
derived similarly as in the previous example; in this context, we
also refer to ~\cite[Remarque 6.4]{cha}. It seems to be an
interesting problem to  consider such kinds of operators in
$L^{p}({\mathbb R}^{n})$ spaces, $p\in [1,\infty),\ p\neq 2,$ by the
use of Fourier multiplier type theorems.

\section{Appendix}

We collect in this appendix results related to $K$-convoluted
$C$-semigroups and cosine functions (\cite{koj}).

\begin{defn} \label{1.-1} \cite{mn}
Let $A$ be a closed operator and $K$ be a locally integrable
function on $[0,\tau )$, $0<\tau \leq \infty $. If there exists a
strongly continuous operator family $(S_K(t))_{t\in [0,\tau)}$ such
that  $S_{K}(t)C=CS_{K}(t),$ $S_{K} (t)A \subset AS_{K}(t),$ $ \int
\limits^{t}_{0} S_{K}(s)xds\in D(A),\; t\in[0,\tau), \;x\in E$ and
\begin{equation} \label{ci}
A\int \limits^{t}_{0} S_{K}(s)xds = S_{K}(t)x-\Theta (t)Cx, \;x\in
E,
\end{equation}
then $(S_K(t))_{t\in [0,\tau)}$ is called a (local) $K$-convoluted
$C$-semigroup having $A$ as a subgenerator. If $\tau=\infty$, then
it is said that $(S_{K}(t))_{t\geq 0}$ is an exponentially bounded,
$K$-convoluted $C$-semigroup with a subgenerator $A$ if, in
addition, there are constants $M>0$ and $\omega \in {\mathbb R}$
such that $||S_{K}(t)||\leq Me^{\omega t},\ t\geq 0.$
\end{defn}

The integral generator of $(S_{K}(t))_{t\in [0,\tau)}$ is defined by
$$
\{(x,y)\in E^{2} : S_{K}(t)x-\Theta(t)Cx=\int
\limits^{t}_{0}S_{K}(s)yds,\ t\in [0,\tau)\}.
$$
It is straightforward to see that the integral generator of
$(S_{K}(t))_{t\in [0,\tau)}$
is an extension of any subgenerator of $(S_{K}(t))_{t\in [0,\tau)}.$

Using Theorem \ref{bau} and ~\cite[Theorem 1.12]{x263} one can prove
the next assertion which can be reformulated for convoluted
$C$-semigroups.

\begin{prop} \label{2.6}
Let $K$ be exponentially bounded and let $A$ be a closed linear
operator with $\{ \lambda^{2} : \mbox{ } Re\lambda
> \omega,\mbox{ }\tilde{K}(\lambda)\neq 0 \} \subset \rho _{C}(A),$
for some $\omega > \max(0,\mbox{abs}(K)).$ Suppose that the function
$$
\lambda \mapsto \lambda
\tilde{K}(\lambda)(\lambda^{2}-A)^{-1}C,\mbox{
}Re\lambda>\omega,\mbox{ }\tilde{K}(\lambda)\neq 0,
$$
can be extended to an analytic function $\Upsilon : \{\lambda:
Re\lambda>\omega\} \rightarrow {\mathbb C}$ satisfying
$$
|\Upsilon (\lambda) | \leq M_0|\lambda|^r , \quad Re\lambda
> \omega,
$$
where $r\geq -1.$ Then, for every $\alpha >1,$ there exist a
continuous function \\ $C_{K}: [0,\infty )\rightarrow L(E)$ with
$C_{K}(0) = 0$ and a constant $M_1>0$ such that\\ $|| C_{K}(t) ||
\leq M_{1}e^{\omega t}$, $ t \geq 0$ and that
$$
\lambda \tilde{K}(\lambda)(\lambda^{2}-A )^{-1}C = \lambda ^{\alpha
+r} \int \limits^{\infty }_{0} e^{-\lambda t} C_{K}(t)dt , \quad
Re\lambda
> \omega,\mbox{   }\tilde{K}(\lambda)\neq 0.
$$
Furthermore, $(C_{K}(t))_{t\geq 0}$ is a norm continuous, $(K\ast
\frac{t^{\alpha+r-1}}{\Gamma(\alpha +r)})$-convoluted $C$-cosine
function with a subgenerator $A.$
\end{prop}

Now we give an improvement of ~\cite[Theorem 3.3]{koj} and
~\cite[Theorem 3.3]{mn}. The similar assertion holds for convoluted
$C$-cosine functions.

\begin{thm}\label{arendt}
Let $K$ satisfy (P1) and let $A$ be a closed linear operator.
\begin{itemize}
\item[(i)] Then $A$ is a subgenerator of an exponentially bounded $\Theta
$-convoluted \\ $C$-semigroup $(S_{\Theta }(t))_{t\geq 0}$
satisfying the condition
\begin{equation}\label{arendt}
||S_{\Theta}(t+h)-S_{\Theta}(t)||\leq Che^{\omega (t+h)},\ t\geq 0 ,
\ h\geq 0,
\end{equation}
$$\mbox{ for some } C>0 \; \mbox{ and }\; \omega \geq 0,$$
if and only if there exists an $a\geq \max (\omega ,
\mbox{abs}(K))$ such that
\begin{equation}\label{arendt1}
\{ \lambda \in {\mathbb C} : Re\lambda >a,\mbox{
}\tilde{K}(\lambda)\neq 0 \} \subset \rho_{C} (A),
\end{equation}
\begin{equation}\label{arendt0}
\lambda \mapsto \tilde{K}(\lambda)(\lambda-A)^{-1}C,\ Re\lambda
>a,\ \tilde{K}(\lambda) \neq 0 \mbox{ is analytic, and}
\end{equation}

\begin{equation}\label{arendt2}
\bigg\| \frac{d^{k}}{d\lambda ^{k}} [\tilde{K}(\lambda
)(\lambda-A)^{-1}C] \bigg\| \leq \frac{Mk!}{(\lambda -\omega
)^{k+1}} , \; k\in {{\mathbb N}_{0}} ,\;  \lambda > a,\mbox{
}\tilde{K}(\lambda)\neq 0.\;
\end{equation}

\item[(ii)] Assume additionally that $A$ is densely defined. Then $A$
is a subgenerator of an exponentially bounded $K$-convoluted
$C$-semigroup $(S_{K}(t))_{t\geq 0}$ satisfying $\| S_{K}(t)\| \leq
Me^{\omega t}, \;t\geq 0,\ \omega \geq 0,$ if and only if there is
an $a\geq \max ( \omega , \mbox{abs}(K) )$ such that
$(\ref{arendt1}),$ $(\ref{arendt0})$ and $(\ref{arendt2})$ are
fulfilled.
\end{itemize}
\end{thm}
{\em Proof:}
 (i) Let us assume $(\ref{arendt1}),$ $(\ref{arendt0})$ and $(\ref{arendt2}).$
Put $a=\max ( \omega , \mbox{abs}(K) ).$ If $\lambda >a\mbox{ and
}\tilde{K}(\lambda)\neq 0 ,$ then $(\ref{arendt2})$ implies that the
power series
$$
\sum_{k\geq 0} \frac{[ \tilde{K}(\lambda )(\lambda-
A)^{-1}C]^{(k)}(\lambda)}{k!}(z-\lambda)^{k},
$$
converges for every $z\in {\mathbb C}$ satisfying
$|z-\lambda|<\lambda -\omega$. It, in turn, implies that there
exists a $C^{\infty}$-function $\Upsilon : (a, \infty) \rightarrow
L(E)$ satisfying $\Upsilon(\lambda)=\tilde{K}(\lambda)(\lambda
-A)^{-1}C,\ \lambda>a,\ \tilde{K}(\lambda)\neq 0$ and
$||\frac{d^{k}}{d\lambda^{k}}\Upsilon(\lambda)||\leq
\frac{Mk!}{(\lambda-\omega)^{k+1}},\ k\in {{\mathbb N}_{0}},\
\lambda>a.$ An application of ~\cite[Corollary 3.3]{hibi} gives that
there exist a constant $C>0$ and a function $S_{\Theta} : [0,\infty)
\rightarrow L(E)$ such that (\ref{arendt}) holds and that
$\Upsilon(\lambda)=\lambda\int \limits^{\infty}_{0}e^{-\lambda
t}S_{\Theta}(t)dt,\ \lambda>a.$ Then it is straightforward to see
that $\tilde{\Theta}(\lambda)(\lambda -A)^{-1}C=\int
\limits^{\infty}_{0}e^{-\lambda t}S_{\Theta}(t)dt,\ Re\lambda>a,\
\tilde{\Theta}(\lambda) \neq 0.$ This implies that
$(S_{\Theta}(t))_{t\geq 0}$ is an exponentially bounded,
$\Theta$-convoluted $C$-semigroup with a subgenerator $A.$ Assume
conversely that $A$ is a subgenerator of an exponentially bounded,
$\Theta$-convoluted $C$-semigroup $(S_{\Theta}(t))_{t\geq 0}$ which
satisfies (\ref{arendt}). Proceeding as before, one obtains
(\ref{arendt1}) and
$$
\lambda (\lambda -A)^{-1}Cx=\frac{1}{\tilde{K}(\lambda )}
\int\limits_{0}^{\infty }e^{-\lambda t}S_{\Theta}(t)xdt,\; x\in E,
\; Re \lambda >a,\mbox{ }\tilde{K}(\lambda)\neq 0.
$$
This implies (\ref{arendt0}). To prove (\ref{arendt2}), let $x\in E$
and $ x^{\ast}\in E^{\ast}$ be fixed. Put now
$f(t):=x^{\ast}(S_{\Theta}(t)x),$ $ t\geq 0.$ Then (\ref{arendt})
implies that $f$ is differentiable almost everywhere in $[0,\infty)$
with $|f^{\prime}(t)|\leq C||x|| ||x^{\ast}||e^{\omega t},\mbox{ for
a.e. } t\geq 0.$ Moreover,
$$
x^{\ast}(\tilde{K}(\lambda)(\lambda-A )^{-1}Cx)=\int
\limits^{\infty}_{0}e^{-\lambda t}f^{\prime}(t)dt,\ \lambda>a,\
\tilde{K}(\lambda)\neq 0.
$$
Therefore, (\ref{arendt2}) is true. Using the same arguments as in
the proof of ~\cite[Theorem 3.4, p. 14]{x263}, one obtains (ii).

The proof of the statements (a), (c) and (d) of the following
theorem is given in \cite{koj} while the proof of (b) follows in
exactly the same way as in the proof of ~\cite[Proposition
1.3]{ksho}.

\begin{prop}\label{natalija}
Suppose $A$ is a subgenerator of a (local) $K$-convoluted
$C$-semigroup $(S(t))_{t\in \lbrack 0,\tau )}$. Let $B$ be the
integral generator of $(S(t))_{t\in \lbrack 0,\tau )}.$ Then:
\begin{itemize}
\item[(a)] $$ S(t)S(s) = \left[ \int\limits^{t+s}_{0}
-\int\limits^{t}_{0} -\int\limits^{s}_{0} \right] K(t+s-r) S(r)C
dr,\;\  0\leq t, \ s, \ t+s < \tau.$$
\item[(b)] $B=C^{-1}BC.$

\noindent If, additionally, $K$ is a kernel, then the next
conditions are satisfied:

\item[(c)] $B=C^{-1}AC.$

\item[(d)] For every $\lambda \in \rho _{C}(A)$ : $(\lambda -A)^{-1}
CS(t) = S(t) (\lambda -A)^{-1}C , \; t\in [0,\tau ) .$
\end{itemize}
\end{prop}


\end{document}